\documentclass[12pt,a4paper]{article}
\usepackage{graphicx}
\usepackage{amsmath}
\usepackage{amssymb}
\usepackage{amsthm}
\usepackage{geometry}
\usepackage{fancyhdr}
\usepackage{color} 
\usepackage{mathabx}
\usepackage{overpic}

\usepackage[all]{xy}
\newcommand{\al}{\alpha}

\newcommand{\be}{\beta}
\newcommand{\ga}{\gamma}

\newcommand{\ra}{\rightarrow}

\newcommand{\xra}{\xrightarrow}

\newcommand{\rgl}{\rangle}
\newcommand{\lgl}{\langle}

\newcommand{\bpf}{\begin{proof}}
\newcommand{\epf}{\end{proof}}
\newcommand{\bthm}{\begin{thm}}
\newcommand{\ethm}{\end{thm}}
\newcommand{\bprop}{\begin{prop}}
\newcommand{\eprop}{\end{prop}}
\newcommand{\bcor}{\begin{cor}}
\newcommand{\ecor}{\end{cor}}
\newcommand{\blem}{\begin{lem}}
\newcommand{\elem}{\end{lem}}
\newcommand{\bdefn}{\begin{defn}}
\newcommand{\edefn}{\end{defn}}
\newcommand{\bexmp}{\begin{exmp}}
\newcommand{\eexmp}{\end{exmp}}
\newcommand{\brem}{\begin{rem}}
\newcommand{\erem}{\end{rem}}

\newcommand{\bdia}{\begin{displaymath}\xymatrix}
\newcommand{\edia}{\end{displaymath}}
\newcommand{\beq}{\begin{equation*}\begin{aligned}}
\newcommand{\eeq}{\end{aligned}\end{equation*}}

\newcommand{\intg}{\mathbb{Z}}
\newcommand{\real}{\mathbb{R}}

\newcommand{\shm}{\underline{\rm SHM}}

\newtheorem{thm}{\textbf {Theorem}}[section]
\newtheorem{cor}[thm]{\textbf{Corollary}}
\newtheorem{prop}[thm]{\textbf{Proposition}}
\newtheorem{lem}[thm]{\textbf{Lemma}}
\newtheorem{conj}[thm]{Conjecture}

\newtheorem{quest}[thm]{Question}

\theoremstyle{definition}
\newtheorem{defn}[thm]{\textbf{Definition}}

\newtheorem{exmp}[thm]{Example}

\theoremstyle{remark}
\newtheorem{rem}[thm]{Remark}

\title{Gluing maps and cobordism maps for sutured monopole Floer homology}

\author{Zhenkun Li}

\date{}

\begin{document}
\bibliographystyle{plain}

\maketitle
\begin{abstract}
The naturality of sutured monopole Floer homology, which was introduced by Kronheimer and Mrowka \cite{kronheimer2010knots}, is an important question and is partially answered by Baldwin and Sivek \cite{baldwin2015naturality}. In this paper we construct the cobordism maps for sutured monopole Floer homology, thus improve its naturality. The construction can be carried out for sutured instantons as well. In the paper we also construct gluing maps in sutured monopoles and sutured instantons.
\end{abstract}

\tableofcontents
\newpage

\section{Introduction}
\subsection{Main theorems and backgrounds}
Sutured manifold is a powerful tool introduced by Gabai \cite{gabai1983foliations} in 1983, to study the topology of $3$-manifolds. In 2010, the construction of monopole Floer homology was carried out on balanced sutured manifold by Kronheimer and Mrowka \cite{kronheimer2010knots}. The combination of Floer theories and sutured manifolds has many important applications. For example, sutured Floer homology can detect tautness (see \cite{juhasz2008floer} and \cite{kronheimer2010knots}) and fibredness of knots (see \cite{ni2007knot} and \cite{kronheimer2010knots}), and these played essential roles in the proofs that Khovanov homology detects unknots by Kronheimer and Mrowka \cite{kronheimer2011khovanov} and that Khovanov homology detects trefoil by Baldwin and Sivek \cite{baldwin2018khovanov}. In this paper, we construct the gluing maps and cobordism maps for sutured monopole and instanton Floer homology. This will enrich our tool bar for potential usage.

\bthm(Gluing maps)
Suppose $(M,\ga)$ and $(M',\ga')$ are balanced sutured manifolds and $M\subset{\rm int}(M)$. Suppose there is a contact structure $\xi$ on $Z=M'\backslash{\rm int}(M)$ so that $\partial{Z}$ is convex with $\ga\cup\ga'$ being the dividing set, then there is a contact gluing map
$$\Phi_{\xi}:\underline{\rm SHM}(-M,-\ga)\ra\underline{\rm SHM}(-M',-\ga'),$$
which is  well defined up to the multiplication by a unit. Furthermore, the gluing map satisfies the following properties:

(1). If $Z\cong \partial{M}\times [0,1]$ then there exists a diffeomorphism
$$\phi:M\ra M',$$
which restricts to the identity outside a collar of $\partial{M}\subset M$ and is isotopic to the inclusion $M\hookrightarrow M'$, so that
$$\Phi_{\xi}\doteq\underline{\rm SHM}(\phi).$$
Here $\doteq$ means equal up to multiplication by a unit.

(2). Suppose $(M'',\ga'')$ is another balanced sutured manifold and $M'\subset{\rm int}(M'')$, and let $Z'=M''\backslash{\rm int}(M')$ with a contact structure $\xi''$ on $Z'$ so that $\partial{Z}'$ is convex with dividing set $\ga'\cup\ga''$, then we have
$$\Phi_{\xi\cup\xi'}\doteq\Phi_{\xi'}\circ\Phi_{\xi}.$$

(3). If $Z\cong \partial{M}\times[0,1]\cup h$ where $\partial{M}\times\{0\}$ is identified with $\partial{M}\subset M$ and $h$ is a contact handle attached to $M\cup Z$ along $\partial{M}\times\{1\}$, then there is a suitable diffeomorphism
$$\phi:M\ra M\cup \partial{M}\times[0,1],$$
which restricts to the identity outside a collar of $\partial{M}\subset M$, and is isotopic to the inclusion $M\hookrightarrow M\cup\partial{M}\times[0,1]$, so that
$$\Phi_{\xi}\doteq C_{h}\circ \underline{\rm SHM}(\phi),$$
where $C_h$ is the contact handle map associated to $h$.
\ethm

\bthm(Sutured cobordism maps)
	Suppose $\mathcal{W}=(W,Z,[\xi])$ is a sutured cobordism between two balanced sutured manifolds $(M_1,\ga_1)$ and $(M_2,\ga_2)$ then $\mathcal{W}$ induces a cobordism map
$$\shm(\mathcal{W}):\underline{\rm SHM}(M_1,\ga_1)\ra\underline{\rm SHM}(M_2,\ga_2),$$
which is well defined up to multiplication by a unit and satisfies the following properties:

(1). Suppose $\mathcal{W}=(M\times[0,1],\partial{M}\times[0,1],[\xi_0])$ so that $\xi_0$ is $[0,1]$-invariant, then
$$\shm(\mathcal{W})\doteq id.$$

(2). Suppose $\mathcal{W}'=(W',Z',[\xi]')$ is another sutured cobordism from $(M_2,\ga_2)$ to $(M_3,\ga_3)$, then we can compose them to get a cobordism $\mathcal{W}''=(W\cup W',Z\cup Z',[\xi\cup \xi'])$ from $(M_1,\ga_1)$ to $(M_3,\ga_3)$ and there is an equality
$$\shm(\mathcal{W}'')\doteq\shm(\mathcal{W}')\circ \shm(\mathcal{W}).$$

(3). For any balanced sutured manifold there is a canonical pairing
$$\lgl\cdot,\cdot\rgl:\underline{SHM}(M,\ga)\times \underline{SHM}(-M,\ga)\ra \mathcal{R},$$
which is well defined up to multiplication by a unit. Here $\mathcal{R}$ is the coefficient ring. Furthermore, let $\mathcal{W}=(W,Z,[\xi])$ be a cobordism from $(M_1,\ga_1)$ to $(M_2,\ga_2)$, and let $\mathcal{W}^{\vee}=(W,Z,[\xi])$ be the cobordism with same data but viewed as from $(-M_2,\ga_2)$ to $(-M_1,\ga_1)$. Then $\mathcal{W}$ and $\mathcal{W}^{\vee}$ induce cobordism maps which are dual to each other under the canonical pairings.
\ethm

A sutured manifold $(M,\ga)$ is a compact oriented $3$-manifold $M$, with an oriented $1$-submanifold $\ga$ on the boundary $\partial{M}$. The $1$-submanifold $\ga$ is called the suture and it divides $\partial{M}$ into two parts $R_+(\ga)$ and $R_-(\ga)$ according to the orientation induced by $\ga$ and $M$. It is called balanced if every component of $\partial{M}$ contains a suture and $R_+(\ga)$ and $R_-(\ga)$ have the same Euler characteristics.

The monopole Floer homology of a closed $3$-manifold $Y$ together with a spin${}^c$ structure $\mathfrak{s}$ on $Y$ was built by a version of infinite dimensional Morse theory based on the so called Chern-Simons-Dirac functional as in \cite{kronheimer2007monopoles}. The total homology group is then a direct sum among all spin${}^c$ structures. To adapt the construction to balanced sutured manifolds, Kronheimer and Mrowka constructed a pair $(Y,R)$, consisting of a closed $3$-manifold and a distinguishing surface, out of the sutured data $(M,\ga)$. This pair was called a closure and the sutured monopole Floer homology was defined to be the monopole Floer homology of $Y$ using only top spin${}^c$ structures with respect to $R$, i.e., those spin${}^c$ structures $\mathfrak{s}$ on $Y$ so that
$$c_1(\mathfrak{s})[R]=2g(R)-2.$$
This homology is denoted by $\shm(M,\ga)$.

The monopole Floer homology of a closed $3$-manifold has very good naturality. This is partially because the space of all choices in the construction of monopole Floer homology is contractible. Hence it is natural to ask whether sutured monopole would also have a good naturality property. However, the construction of the closure, which involves some 'discrete' choices, make the question much more difficult to be studied. The naturality of sutured monopole is partially proved by Baldwin and Sivek in \cite{baldwin2015naturality}, where they showed that for a fixed balanced sutured manifold and any two different closures of it, there is a canonical map for Floer homologies between them. However, whether there exists a cobordism map in sutured monopole Floer homology theory is still open (and the main theorem of this paper answers this question positively). To be compared with, Juh\'asz \cite{juhasz2016cobordisms} has constructed a cobordism map in sutured (Heegaard) Floer homology theory.

Our construction of the cobordism map would provide the sutured monopole Floer homology a better naturality property. In some cases when we could fix the choices of closures, we might be able to make use of the even better naturality of monopole Floer homology of closed $3$-manifolds and find some future applications.

Along with the cobordism map, we also construct gluing maps in sutured monopoles and instantons. Gluing maps are very important tools in the sutured Floer homology theory. One direct application of gluing maps in future is to construct a possible minus version of knot Floer homology in monopole or instanton settings, using a direct system whose morphisms coming from gluing maps. The question has been proposed in details in the introductory part of Baldwin and Sivek's paper \cite{baldwin2016contact}. One thing we would like to comment here is that the construction of the direct system would be an immediate application of the gluing maps, but the construction of gradings in the direct limit would have some difficulties, and this will be the main topic of the author's future paper \cite{li2019direct}.

\subsection{Outline of the proof}
The topological data for a cobordism between two balanced sutured manifolds $(M_1,\ga_1)$ and $(M_2,\ga_2)$ would be a pair $(W,Z)$ where $W$ is a compact oriented $4$-manifold with boundary
$$\partial{W}=-M_1\cup Z\cup M_2.$$
In order to keep track of the sutured data, we need also a contact structure $\xi$ on $Z$ so that $\partial{Z}$ is convex and $\ga_1\cup\ga_2$ is the dividing set. For the purpose of gluing cobordisms, we should allow the contact structure $\xi$ to vary by isotopy and look at only the isotopy class $[\xi]$ of contact structures on $Z$. We call $\mathcal{W}=(W,Z,[\xi])$ a sutured cobordism from $(M_1,\ga_1)$ to $(M_2,\ga_2)$.

The construction of the cobordism map would share some similarity of the construction of that in sutured Heegaard Floer homology theory in Juh\'asz \cite{juhasz2016cobordisms}. The construction falls into two steps and the first is to use the isotopy class of contact structures $[\xi]$ on $Z$ to construct a gluing map
$$\Phi_{-\xi}:\underline{\rm SHM}(M_1,\ga_1)\ra \underline{\rm SHM}(M_1\cup(-Z),\ga_2).$$
A corresponding construction for sutured (Heegaard) Floer homology was done by Honda, Kazez and Mati\'c \cite{honda2008contact} and revisited by Juh\'asz and Zemke \cite{juhasz1803contact}. The second step is to construct a map from the $4$-manifold $W$
$$F_W: \underline{SHM}(M_1\cup(-Z),\ga_2)\ra \underline{SHM}(M_2,\ga_2).$$
The composition $\shm(\mathcal{W})=F_W\circ\Phi_{-\xi}$ would be the desired cobordism map.

The second map $F_{W}$ arising from $W$ is straightforward in monopole settings. The main difficulty is to construct the gluing map $\Phi_{\xi}$. There is no known construction prior to this paper. Some partial works were done by Baldwin and Sivek in \cite{baldwin2016contact}, where they only constructed the handle attaching maps for contact handle attachments. A contact handle is a tight contact $3$-ball attached to a balanced sutured manifold and is attached in different ways according to the index of the handle. The straightforward idea then is to decompose $Z$ into contact handles and composite the gluing maps defined for those handles. However they only conjectured that two different contact handle decompositions will result in the same composition map.

In the present paper we are going to introduce slightly different definitions for contact handle attaching maps for $2$- and $3$-handles. Though they turn out to be equivalent to what have been constructed by Baldwin and Sivek \cite{baldwin2016contact}, our point of view will be a little bit more convenient when studying the duality of the sutured cobordism map. In the paper we also make use of a tool called contact cell decomposition, which was introduced by Juh\'asz and Zemke \cite{juhasz1803contact}. A contact cell decomposition can be thought of a refinement of the construction of Legendrian graphs inside contact $3$-manifolds, by Honda, Kazez and Mati\'c \cite{honda2009contact}, as a preparation for defining the contact elements in sutured (Heegaard) Floer homology. In a contact cell decomposition of $Z$, we decompose $Z$ into three pieces $Z=N\cup Z'\cup N'$. Here $N\cup N'$ is a collar of the boundary $\partial{Z}\subset Z$ and $Z'$ is decomposed further into contact handles so that any two different decompositions are related by isotopies and three types of handle cancelations. In this paper we are able to prove that the composition of gluing map is independent of all three types of handle cancelations and hence get a well defined gluing map.

As an application of the gluing map, we prove the following result, which is originally conjectured by Baldwin and Sivek \cite{baldwin2016contact}.

\bcor
Under the above settings, suppose there are two different ways of contact handle decompositions of $Z$, both relative to $\partial{M}$:
$$M'=M\cup h_1\cup...\cup h_n,~M'=M\cup h_1'\cup...\cup h_m'.$$
Then the compositions of the two sets of handle attaching maps are the same:
$$C_{h_n}\circ...\circ C_{h_1}\doteq C_{h_m'}\circ...\circ C_{h_1'}:\shm(-M,-\ga)\ra\shm(-M',-\ga').$$
\ecor

Further with the discussion in \cite{baldwin2016contact}, we know the follow thing.
\bcor
With the above notations, if the contact structure $\xi$ on $Z$ is the restriction of a contact structure $\xi'$ on all of $M'$, then the gluing map preserves contact elements, i.e.,
$$\Phi_{\xi}(\phi_{\xi'|_M})\doteq \phi_{\xi'}.$$
\ecor

When composing the gluing map $\Phi_{-\xi}$ with the map $F_{W}$ coming from the $4$-manifold $W$ as discussed above, we get the cobordism map $\shm(\mathcal{W})$ associated to the sutured cobordism $\mathcal{W}=(W,Z,[\xi])$.

The functoriality of the cobordism maps holds. This is essentially because we can interpret the two maps $\Phi_{-\xi}$ and $F_{W}$ as attaching $4$-dimensional handles to a suitable product cobordism, and we can somehow prove that  handles attached corresponding to different steps can change the order of attaching with each other.

The duality of the cobordism map is also proved. To do this We actually introduced a second way to construct the gluing maps as well as the cobordism maps, so that the duality is then a simple corollary.   

Although we will work with local coefficients though out the paper, we shall remark here that all discussions can be modified to work with simply $\intg$ coefficients (and we shall fix a large enough genus for closures of balanced sutured manifolds) except for proposition \ref{prop_gluing_disjoint_union}. When using $\intg$ coefficients, the ambiguity appeared in the above statements reduces to being up to a sign. The reason why \ref{prop_gluing_disjoint_union} relies on local coefficients is that we shall use Floer excisions along tori, which was introduced in \cite{kronheimer2010knots}, in the proof of that proposition. However local coefficients are necessary in that setting. 

The sutured instanton Floer homology was also introduced by Kronheimer and Mrowka \cite{kronheimer2010knots}. A parallel construction for sutured instanton can also be done in a similar way. We will briefly discuss about sutured instantons in the last section of the paper. It worth mentioning here that in \cite{baldwin2016instanton} Baldwin and Sivek defined the contact elements as well as contact handle gluing maps for sutured instanton. However, they only proved that the contact element is preserved under $0,1,2$-handle attaching maps but didn't say anything about $3$-handles. Using gluing maps constructed in this paper, we are able to prove that contact elements are also preserved by $3$-handle attaching maps. 

\subsection{Future questions}
The construction of gluing maps and cobordism maps would be a first step to many further problems and we would like to introduce some of them here. We have already mentioned one above on the minus version of knot monopole Floer homology and here are more questions to be asked.

 A first adaption of the construction in the current paper might be to the sutured knot (or link) homology. As suggested by Juh\'asz \cite{juhasz2016cobordisms}, we are given a cobordism $(X,F,\sigma)$ between two links $L_1\subset Y_1$ and $L_2\subset Y_2$ with marked points, where $X$ is a cobordism between $Y_1$ and $Y_2$, $F\subset X$ is a cobordism between $L_1$ and $L_2$ and $\sigma\subset F$ is a $1$-dimensional submanifold determining a contact structure on the boundary of a tubular neighborhood of $F\subset X$. Then we could try to construct a map between link Floer homologies. There may be some further applications to the study of embedded surfaces in $4$-manifolds.

A second question is related to contact elements of the sutured manifolds. Given a balanced sutured manifold $(M,\ga)$ with a contact structure $\xi$ so that $\partial{M}$ is convex and $\ga$ is the dividing set, Baldwin and Sivek constructed in \cite{baldwin2016contact} a closure $(Y,R)$ of $(M,\ga)$ which carries a contact structure $\bar{\xi}$ restricting to $\xi$ on $M\backslash N(\ga)$. Here $N(\ga)$ is a neighborhood of $\ga\subset M$. Hence they were able to define a contact element
$$\phi_{\xi}\in\underline{SHM}(-M,-\ga)$$
based on work by Kronheimer, Mrowka, Ozsv\'ath and Szab\'o \cite{kronheimer2007monopolesandlens}. However they only carried out this construction using connected auxiliary surface (an auxiliary surface is the surface used to construct closures of balanced sutured manifolds), while in some cases, disconnected auxiliary surfaces might be more convenient (see \cite{kronheimer2010knots}, section 6.) So it would be interesting to generalize their construction using disconnected auxiliary surface and study how contact invariants behave under the Floer excision maps defined in \cite{kronheimer2010knots}. Another related questions is that in \cite{baldwin2016contact}, or in second 6 of the current paper, contact invariants for sutured instantons are defined. One can ask what is the analytical correspondence in the classical Instanton Floer homology theory.

A third question is about trace and co-trace cobordisms. Suppose $(M,\ga)$ is a balanced sutured manifold, then we can form a special cobordism 
$$\mathcal{W}=(M\times[0,1],\partial{M}\times[0,1],[\xi_0])$$
where $\xi_0$ is a contact structure on $\partial{M}\times[0,1]$ so that $\xi_0$ is $[0,1]$-invariant, $\partial{M}\times\{t\}$ is convex for each $t\in[0,1]$ and $\ga\times\{t\}$ is the dividing set. We can view $\mathcal{W}$ as a cobordism from $(M\sqcup(-M),\ga\cup\ga)$ to $\emptyset$. In \cite{juhasz1803contact} the corresponding cobordism map for sutured (Heegaard) Floer homology was computed and one could ask whether we have a similar result for sutured monopoles.

A forth question is about the ambiguity of being up to multiplication by a unit. When using $\intg$ coefficients, it is up to a sign, which is kind of acceptable, as the contact invariant is also only defined up to a sign and this ambiguity cannot be resolved, as shown by Lin \cite{lin2018froyshov}. Also, as a comparison, Honda, Kazez and Matic's construction of gluing map in sutured (Heegaard) Floer homology also has a sign ambiguity when using $\intg$ coefficients. However, when using general local coefficient over a suitable ring $\mathcal{R}$, it might not be satisfactory. For example, the pairing
$$\lgl \cdot,\cdot \rgl: \underline{\rm SHM}(M,\ga)\times \underline{\rm SHM}(-M,\ga)\ra \mathcal{R}$$
defined above is also up to a unit. In the worst case where $\mathcal{R}$ is a field, we only know that the vanishing or non-vanishing of the pairing is well defined. So it would be interesting to see whether or not one could improve this ambiguity.

The paper is organized as follows. In section 2 we review the basic settings of the sutured monopole Floer homology and the naturality. In section 3, we discuss on the construction of contact handle attaching maps and prove their cancelation or invariance properties. Those basic ingredients then are used in section 4 for constructing general gluing maps and proving basic properties of them. In section 5, we construct the cobordism maps associated to sutured cobordisms between balanced sutured manifolds and prove their basic properties. In section 6 we briefly go through the construction in sutured instanton Floer homology.

{\bf Acknowledgements}. This material is based upon work supported by the National Science Foundation under Grant No. 1808794. I would like to express my enormous gratitude towards my advisor Tomasz Mrowka for suggesting the problem and the invaluable helps all the way along. I would like to thank John Baldwin, Mariano Echeverria, Jianfeng Lin, Yu Pan and Boyu Zhang for helpful conversations.

\section{Prelimilaries}
\subsection{Monopole Floer homology for $3-$manifold}
Suppose $(Y,\mathfrak{s})$ is a closed connected oriented 3-manifold equipped with a spin${}^c$ structure $\mathfrak{s}$. Kronheimer and Mrowka in \cite{kronheimer2007monopoles} associated 3 flavors of monopole Floer homologies to $(Y,\mathfrak{s})$, with $\intg$ coefficients:
$$\widehat{HM}_{\bullet}(Y,\mathfrak{s}),~\widecheck{HM}_{\bullet}(Y,\mathfrak{s}),~\widebar{HM}_{\bullet}(Y, \mathfrak{s}).$$
The three flavors fit into a long exact sequence:
\begin{equation}\label{eq_long_exact_sequence}
...\ra\widebar{HM}_{\bullet}(Y, \mathfrak{s})\xra{i} \widecheck{HM}_{\bullet}(Y, \mathfrak{s})\xra{j} \widehat{HM}_{\bullet}(Y, \mathfrak{s})\xra{p} \widebar{HM}_{\bullet}(Y, \mathfrak{s})\ra...
\end{equation}

Suppose we are given a smooth $1$-cycle $\eta\subset Y$, and let 
$\mathcal{R}$ be the Novikov ring over $\intg$, which is defined as
$$\mathcal{R}=\{\sum_{\al}n_{\al}t^{\al}|\al\in\real,~n_{\al}\in\intg,~\sharp\{\al\in\real|n_{\al}<N\}<\infty~for~all~N\in \intg\}.$$
Then we can define similarly all three flavors of monopole Floer homologies with local coefficients
$$\widehat{HM}_{\bullet}(Y,\mathfrak{s};\Gamma_{\eta}),~\widecheck{HM}_{\bullet}(Y,\mathfrak{s};\Gamma_{\eta}),~\widebar{HM}_{\bullet}(Y, \mathfrak{s};\Gamma_{\eta}).$$
They also fit into a the same long exact sequence as (\ref{eq_long_exact_sequence}). 

If furthermore the spin${}^c$ structure is non-torsion, that is, $c_1(\mathfrak{s})$ is not a torsion element in $H^2(M;\intg)$, then we have
$$\widebar{HM} _{\bullet}(Y,\mathfrak{s};\Gamma_{\eta})=0,$$
and $\widecheck{HM} _{\bullet}(Y,\mathfrak{s};\Gamma_{\eta})$ is isomorphic to $\widehat{HM} _{\bullet}(Y,\mathfrak{s};\Gamma_{\eta})$ via $j$. So we will call either flavor to be just $HM_{\bullet}(Y,\mathfrak{s};\Gamma_{\eta})$.

Suppose $F\subset Y$ is a closed oriented embedded surface of genus at least 2. Let $\mathfrak{S}(Y|F)$ be the set of all spin${}^c$ structures $\mathfrak{s}$ such that
$$c_1(\mathfrak{s})[F]=2g(F)-2,$$
and define
$$HM(Y|F;\Gamma_{\eta})=\bigoplus_{\mathfrak{s}\in\mathfrak{S}(Y|F)}HM _{\bullet}(Y,\mathfrak{s};\Gamma_{\eta}).$$

Suppose $(Y_1,F_1,\eta_1)$ and $(Y_2,F_2,\eta_2)$ are two triples, then a cobordism $(W,F_W,\nu)$ between them is a triple where

(1). $W$ is a cobordism from $Y_1$ and $Y_2$, which means that $W$ is a smooth compact oriented $4$-manifold with boundary and there is an orientation preserving diffeomorphism from $\partial{W}$ to $-Y_1\cup Y_2$.

(2). $F_W\subset W$ is a closed oriented embedded surface in $W$, which contains $F_1$ and $F_2$ as two components. 

(3). We have $\nu\subset W$ being a smooth $2$-cycle and $\partial{\nu}=(-\eta_1)\cup\eta_2$

As discussed in \cite{kronheimer2010knots}, the cobordism $(W,F_W,\nu)$ induces a map between mononopole Floer homologies:
$$HM(W|F_W;\Gamma_{\nu}):HM(Y_1|F_1;\Gamma_{\eta_1})\ra HM(Y_2|F_2;\Gamma_{\eta_2}).$$

\brem
For simplicity, in the rest of the paper, we may omit the surface and local coefficients from the notation of a cobordism map.
\erem

\subsection{Sutured monopole Floer homology}
\begin{defn}
A {\it balanced sutured manifold} $(M,\ga)$ consists of the following data:

(1). A compact, oriented 3-manifold $M$ with non-empty boundary $\partial{M}$.

(2). An embedded oriented 1-submanifold $\ga\subset \partial{M}$. 

(3). An annular neighborhood $A(\ga)$ of $\ga$ on $\partial{M}$, which can be identified with $\ga\times [-1,1]$.

(4). $R(\ga)=\partial{M}\backslash\mathring{A}(\ga)$ being the closure of the complement of $A(\ga)$ on $\partial{M}$.

They should satisfy the following requirements:

(1). $M$ has no closed components.

(2). Every component of $\partial{M}$ contains at least one suture.

(3). $R(\ga)$ can be oriented in a way that $\partial{R(\ga)}$, as oriented curve, is parallel to $\gamma$ in $A(\ga)$. The requirement (2) above makes sure that $R(\ga)$ has no closed components so the orientation above is unique, and is called the {\it canonical orientation}. 

(4). Let $R_+(\ga)$ be the part of $R(\ga)$ so that the canonical orientation induced by $\ga$ coincides with the boundary orientation of $M$, and let $R_-(\ga)=R(\ga)\backslash R_+(\ga)$. Se shall require further that
$$\chi(R_+(\ga))=\chi(R_-(\ga)).$$
\end{defn}

Now suppose $(M,\ga)$ is a balanced sutured manifold. In order to define the sutured monopole Floer homology, we need to construct a closed 3-manifold with a distinguishing surface inside it. To do this, we pick a compact oriented surface $T$ so that

(1). We have $g(T)\geq 2$.

(2). There exists an orientation reversing diffeomorphism
$$f:\partial{T}\ra \ga.$$

(3). There is a curve $c\subset T$ so that $c$ represents a non-trivial class in $H_1(T)$.

Since $A(\ga)$ has been identified with $\ga\times [-1,1]$, we have a map
$$f\times id: \partial{T}\times [-1,1]\ra \ga\times [-1,1]=A(\ga).$$
We can use this map to glue $T\times [-1,1]$ to $M$:
$$\widetilde{M}=M\mathop{\cup}_{f\times id}T\times [-1,1].$$

The boundary of $\widetilde{M}$ consists of two components
$$R_+=R_+(\ga)\cup(T\times \{1\}),$$
$$R_-=R_-(\ga)\cup(T\times \{-1\}).$$ 

Let $h:R_{+}\ra R_{-}$ be an orientation preserving diffeomorphism so that
$$h(c\times\{1\})=c\times\{-1\}.$$
We can use $id$ and $h$ to glue $R_+\times [-1,1]$ to $\widetilde{M}$ to get a closed manifold $Y$. In details, $R_+\times \{-1\}$ is glued to $R_+\subset \partial{\widetilde{M}}$ via identity and $R_{+}\times \{1\}$ is glued to $R_-\subset\partial{\widetilde{M}}$ via $h$. Let $R$ be the surface $R_+\times \{0\}$. There is then a curve $c=c\times\{0\}\subset R$.

Based on this construction, we have the following definition.

\bdefn\label{defn_closure}
In the above construction, we call $T$ an {\it auxiliary surface} and $h$ a {\it gluing diffeomorphism}. We call the manifold $\widetilde{M}$ a {\it pre-closure} of $(M,\ga)$, and call the pair $(Y,R)$ a {\it closure}. We define the {\it genus} of the closure $(Y,R)$ to be the genus $g(R)$ of $R$.
\edefn

\brem
The definition of genus follows from Baldwin and Sivek \cite{baldwin2015naturality}. While the others follow from Kronheimer and Mrowka \cite{kronheimer2010knots}
\erem

\bdefn\label{defn_sutured_homology}
Suppose $(Y,R)$ is a closure of a balanced sutured manifold $(M,\ga)$ and $\eta$ is the non-separating oriented smooth curve defined as above. Then we define the {\it sutured monopole Floer homology with local coefficients} of $(M,\ga)$ to be
$$SHM(M,\ga;\Gamma_{\eta})=HM(Y|R;\Gamma_{\eta}).$$
\edefn

\subsection{The naturality of sutured monopole Floer homology}
In the definition of sutured monopole Floer homology, there are a few choices $(T,f,h)$ involved (also $c$ and $\eta$). In \cite{kronheimer2010knots} Kronheimer and Mrowka have already proved the invariance:
\bthm
The isomorphism class of sutured monopole Floer homology of a fixed sutured manifold $(M,\ga)$ is independent of all the choices made in the construction of the closure as in definition \ref{defn_sutured_homology}.
\ethm

Although we have the invariance of the isomorphism types of the sutured monopole Floer homologies, it is still not enough to talk about elements in them. This leads to Baldwin and Sivek's work on the naturality of sutured monopole Floer homologies. To get the naturality, Baldwin and Sivek defined a more refined version of closures in \cite{baldwin2015naturality}:

\bdefn\label{defn_marked_closure}
Suppose $(M,\ga)$ is a balanced sutured manifold then a {\it marked closure} of $(M,\ga)$ is a quintuple $\mathcal{D}=(Y,R,r,m,\eta)$ where

(1). $Y$ is a closed oriented smooth $3-$manifold.

(2). $R$ is a connected closed oriented smooth surface with genus at least $2$.

(3). We have a map
$$r: R\times [-1,1]\ra Y$$
which is a smooth orientation preserving embedding.

(4). We have a map
$$m: M\ra Y\backslash {\rm int}({\rm im}(r))$$
which is a smooth orientation preserving embedding and satisfies following properties:

(a). We have that $m$ extends to a diffeomorphism
$$m: M\cup_{f} T\times [-1,1]\ra Y\backslash {\rm int}({\rm im}(r))$$
for some $A(\ga),T,f$ as defined in definition \ref{defn_closure}. Also we need some smooth structure on $M\cup T\times [-1,1] $ which restricts to the given one on $M$.

(b). We have that $m$ restricts to an orientation preserving embedding
$$m: R_{+}(\ga)\backslash A(\ga)\hookrightarrow r(R\times\{-1\}).$$
(This is to make sure that $R$ has the correct orientation.)

(5). We have that $\eta$ is a non-separating smooth oriented curve on $R$.

We define the {\it genus} of $\mathcal{D}$, which is denoted by $g(\mathcal{D})$, to be the genus of the surface $R$. 

We define the sutured Floer homology of the marked closure $\mathcal{D}$ to be
$$SHM(\mathcal{D})=HM(Y|R;\Gamma_{\eta}).$$
\edefn

\brem
Here, strictly speaking, the surface should be $r(R\times\{0\})$, but for simplicity we will always write $R$ for short. Also the local coefficient should be $\Gamma_{r(\eta\times\{0\})}$ and we will only write $\Gamma_{\eta}$.
\erem

\brem
We shall emphasis here that in the requirement (a), the auxiliary surface $T$ should be connected. This is implicitly contained in Baldwin and Sivek's original construction in \cite{baldwin2016contact}. Especially, when they constructed the handle gluing maps, they used the fact that auxiliary surfaces they used are all connected. So throughout the present paper, when we use an auxiliary surface to construct a closure, it should be understood to be connected otherwise stated.
\erem

Baldwin and Sivek also constructed in \cite{baldwin2015naturality} canonical isomorphisms between the homologies of two different marked closures of a fixed balance sutured manifold. The basic terms are canonical maps $\Phi^{g}_{\mathcal{D},\mathcal{D}'}$ for $g(\mathcal{D})=g(\mathcal{D}')$ and $\Phi^{g,g+1}_{\mathcal{D},\mathcal{D}'}$ for $g(\mathcal{D}')=g(\mathcal{D})+1$. These tow basic types of canonical maps will composite to get canonical maps between any two marked closures of the same balanced sutured manifold. In summary they satisfy the following proposition:

\bprop
Suppose $(M,\ga)$ is a balanced sutured manifold. Then for any two marked closures $\mathcal{D}$ and $\mathcal{D}'$ of $(M,\ga)$, there is a canonical map
$$\Phi_{\mathcal{D},\mathcal{D}'}:SHM(\mathcal{D})\ra SHM(\mathcal{D}'),$$
which is well defined up to multiplication by a unit of $\mathcal{R}$, such that 

(1). If $\mathcal{D}=\mathcal{D}'$, then
$$\Phi_{\mathcal{D},\mathcal{D}'}\doteq id.$$
Here $\doteq$ means equal up to multiplication by a unit.

(2). If there are 3 marked closures $\mathcal{D}$, $\mathcal{D}'$ and $\mathcal{D}''$, then we have
$$\Phi_{\mathcal{D}',\mathcal{D}''}\circ \Phi_{\mathcal{D},\mathcal{D}'}\doteq\Phi_{\mathcal{D},\mathcal{D}''}.$$
\eprop

Hence the homologies and canonical maps fit into what is called a projective transitive system in \cite{baldwin2015naturality}:

\bdefn
A {\it projective transitive system of $\mathcal{R}$-modules} consists of an index set $A$ together with

(1). A collection of $\mathcal{R}$-modules $\{M_{\al}\}_{\al\in A}$

(2). A collection of equivalent classes of $\mathcal{R}$-modules homomorphisms $\{ [h_{\al,\be}]\}_{\al,\be\in A}$, such that

(a). Two morphisms are called equivalent if they differed by multiplication by a unit.

(a). For all $\al,\be\in A$, $h_{\al,\be}$ is an isomorphism from $M_{\al}$ to $M_{\be}$. 

(b). If $\al=\be$, then $h_{\al,\be}\doteq id$.

(c). For all $\al,\be,\ga\in A$, we have
$$h_{\be,\ga}\circ h_{\al,\be}\doteq h_{\al,\ga}.$$  
\edefn

With a projective transitive system, we can construct a canonical projective module out of it:
\bdefn
Suppose $(A,\{M_{\al}\}, \{h_{\al,\be}\})$ is a projective transitive system, then we can define a {\it canonical projective module} or simply a {\it canonical module}:
$$M=\coprod_{\al\in A}M_{\al}\slash\sim,$$
where if we have $m_{\al}\in M_{\al}$ and $m_{\be}\in M_{\be}$, then $m_{\al}\sim m_{be}$ if and only if 
$$\pm h_{\al,\be}(m_{\al})=u\cdot m_{\be}.$$
Here $u\in\mathcal{R}^{\times}$ is a unit.
\edefn

\brem
Although we call $M$ simply a canonical module, it shall be understand that it is not a real module. Note all $M_{\al}$ are isomorphic, and we can regard $M$ as having a bijection to
$$M_{\al}\slash \mathcal{R}^{\times}.$$
for all $\al\in A$.
\erem

We can define the maps between the two systems:
\bdefn
Suppose we have two projective transitive systems $(A,\{M_{\al}\},\{h_{\al,\be}\})$ and $(A',\{M'_{\ga}\},\{h'_{\ga,\delta}\})$. A {\it morphism} between them is a collection of equivalent classes of maps $\{[f_{\al,\ga}]\}_{\al\in A,\ga\in A'}$, where two maps are called equivalent if and only if they differ by multiplication by a unit, such that
$$f_{\be,\ga}\circ h_{\al,\be}\doteq h'_{\ga,\delta}\circ f_{\al,\ga}.$$

Such a morphism will define a map between the canonical projective modules
$$f:M\ra M'$$
by choose any $\al\in A$, $\delta\in A'$ and define
$$f([m_{\al}])=[f_{\al,\delta}(m_{\al})].$$
\edefn

\brem
Strictly speaking, a morphism is a collection of maps but for simplicity, we will write one map in the collection to represent it.
\erem

There is a simple lemma about how to compare such two morphisms:
\blem\label{lem_compare_canonical_maps}
Suppose $\{f_{\al,\ga}\}$ and $\{f'_{\al,\ga}\}$ are two morphisms between projective transitive systems $(A,\{M_{\al}\},\{h_{\al,\be}\})$ and $(A',\{M'_{\ga}\},\{h'_{\ga,\delta}\})$, then the following 3 conditions are equivalent:

(1). The induced maps are equal:
$$f=f': M\ra M'.$$

(2). There exists $\al,\be\in A$ and $\ga\in A'$ so that
$$f_{\be,\ga}\circ h_{\al,\be}\doteq f'_{\al,\ga}.$$

(3). There exists $\al\in A$ and $\ga,\delta\in A'$ so that
$$f_{\al,\delta}\doteq h'_{\ga\delta}\circ f'_{\al,\ga}.$$
\elem

From the above discussion, we know that the marked closures $\{SHM(\mathcal{D})\}$ and the canonical maps $\{[\Phi_{\mathcal{D},\mathcal{D}'}]\}$ together form a projective transitive system and hence we have a canonical projective module
$$\underline{\rm SHM}(M,\ga)$$
associated to it. There is a sub-system of it, namely the system consists of $\{SHM(\mathcal{D})\}$ with $\mathcal{D}$ having a fixed genus $g$, and $\{[\Phi^g_{\mathcal{D},\mathcal{D}'}]\}$ which are canonical maps between marked closures of the same genus. This sub-system can be associated to a canonical projective module
$$\underline{\rm SHM}^g(M,\ga)$$

\brem
Throughout the paper, we will use $SHM$ to denote the homology of a particular marked closure $\mathcal{D}$ or some time, when we only care about the isomorphism class, the homology of a balanced sutured manifold $(M,\ga)$. The notation $\underline{\rm SHM}$ will only be used to denote the canonical module coming from projective transitive system over $(M,\ga)$. This usage might be slightly different from Baldwin and Sivek's original paper.
\erem

For later references, we shall present the definition of the canonical map $\Phi^g_{\mathcal{D},\mathcal{D}'}$ between two marked closures of the same genus. We only introduce the definition here. The well definedness and other basic properties were proved in \cite{baldwin2015naturality}.

\blem\label{lem_unique_isotopy_class}(Baldwin, Sivek, \cite{baldwin2015naturality})
If $\Sigma$ is a closed orientable surface of genus at least two, then the space of all diffeomorphisms from $\Sigma\times [0,1]$ to itself, which restrict to identity on $\Sigma\times\{0,1\}$, is connected.
\elem

Suppose $(M,\ga)$ is a balanced sutured manifold and 
$$\mathcal{D}=(Y,R,r,m,\eta),~\mathcal{D}'=(Y',R',r',m',\eta')$$
are two marked closures of $(M,\ga)$ with the same genus $g(\mathcal{D})=g(\mathcal{D}')$. Pick a diffeomorphism
$$C:Y\backslash{\rm int}({\rm im}(r))\ra Y'\backslash{\rm int}({\rm im}(r')),$$
so that
$$C|_{m(M)}=m'\circ m^{-1}:m(M)\ra m'(M).$$
Define
$$\varphi^C_{\pm}=(r'(\pm1,\cdot))^{-1}\circ C\circ (r(\pm1,\cdot)):R\ra R',$$
and 
$$\varphi^C=(\varphi^C_+)^{-1}\circ (\varphi^C_-): R\ra R.$$
Pick a diffeomorphism $\psi^C:R\ra R$ so that
$$\varphi_-^C\circ\psi^C(\eta)=\eta'.$$

If we pick $t<0<t'$, and cut $Y$ open along $r(R\times\{t\})$ and $r(R\times\{t'\})$, and re-glue by the maps
$$r\circ ((\psi^C)^{-1}\times id_t)\circ r^{-1},~r\circ ((\varphi^C\circ\psi^C)\times id_{t'})\circ r^{-1}$$
respectively, we will get another marked closure $(Y'',R,r'',m,\eta)$. The way we construct $\varphi^C$ and $\psi^C$ will ensure that  $Y''$ and $Y'$ are actually diffeomorphic.

To proceed, we want to construct a cobordism from $Y$ to $Y''$ to define the canonical map. The idea is that we can decompose the two gluing maps as compositions of $\pm1$ Dehn twists and such Dehn twists are related to $\mp1$ Dehn surgeries along the curves on which we perform Dehn twists (See \cite{baldwin2015naturality}, section 4.1). Furthermore, those Dehn surgeries are related to attaching $4$-dimensional handles to a suitable product $4$-manifold. (See Rolfson \cite{rolfsen2003knot} Chapter 9) The resulting $4$-manifold can be viewed as a cobordism between two closures $Y$ and $Y''$, and it thus leads to the canonical map. 

Now suppose $(\psi^{C})^{-1}$ and $\varphi^C\circ\psi^C$ are isotopic to the compositions of Dehn twists:
$$(\varphi^C\circ\psi^C)\sim D_{a_1}^{e_1}\circ...\circ D_{a_n}^{e_n},$$
$$(\varphi^C)^{-1}\sim D^{e_{n+1}}_{a_{n+1}}\circ...\circ D_{a_{m}}^{e_m}.$$
Here $D_{a_{i}}^{e_i}$ means doing a Dehn twist along curves $a_i\subset R$. The sup-script $e_i$ is chosen from $\{-1,1\}$ and $1$ represents a positive Dehn twist (or right handed Dehn twist, see section \cite{farb2011primer} 3.1) while $-1$ represents a negative one.

Pick
\begin{equation}\label{eq_choose_t_i}
-\frac{3}{4}<t_m<...<t_{n+1}<-\frac{1}{4}<\frac{1}{4}<t_n<...<t_1<\frac{3}{4}	,
\end{equation}
and pick $t_i'$ to be greater than $t_i$ and smaller than the next number in the sequence (\ref{eq_choose_t_i}). Define
$$\mathcal{N}=\{i|e_{i}=-1\},~\mathcal{P}=\{i|e_i=1\}.$$

Now let $Y_-$ be the $3$-manifold gotten from $Y$ by doing $(+1)$-surgeries along curves $a_i\times\{t_i\}\subset r(R\times\{t_i\})$ for all the indices $i\in\mathcal{N}$. Then $Y$ is diffeomorphic to the manifold gotten from $Y_-$ by doing $-1$-surgeries along curves $a_i\times\{t_i'\}\subset R\times\{t_i'\}$ for all $i\in \mathcal{N}$. If we require such diffeomorphism to restrict to identity on $Y\backslash(R\times((-\frac{3}{4},-\frac{1}{4})\cup(\frac{1}{4},\frac{3}{4})))$ then by lemma \ref{lem_unique_isotopy_class}, there is a unique isotopy class of such diffeomorphisms. If we attach $-1$-framed $4$-dimensional $2$-handles to $Y_-\times [0,1]$ along curves $a_i\times\{t_i'\}\times\{1\}$ for all $i\in\mathcal{N}$ and let the resulting $4$-manifold be $X_-$, then $X_-$ is a cobordism from $Y_-$ to $Y$. Choose the surface $F_{X_-}$ to be $r(R\times\{0\})\times\{0\}$ and choose the $2$-cycle to be $\nu=r(\eta\times\{0\})\times[0,1]$, we can define a map between monopole Floer homologies:
$$HM(X_-): HM(Y_-|r(R\times\{0\});\Gamma_{\eta})\ra HM(Y|r(R\times\{0\});\Gamma_{\eta}).$$

\brem
The $(+1)$-surgeries above means $+1$ with respect to the surface framing $r(R\times\{t_i\})$. In the rest of the paper, when we do surgery with respect to a surface framing and the surface is understood, we may not mention the choice of framings anymore.
\erem

Now let $Y_+$ be the $3$-manifold obtained from $Y_-$ by doing $(-1)$-surgeries along curves $r(a_{i}\times\{t_i\})$ for all $i\in\mathcal{P}$. Similarly as above, there is a cobordism $X_+$ from $Y_-$ to $Y_+$ and a map
$$HM(X_+): HM(Y_-|r(R\times\{0\}))\ra HM(Y_+|r(R\times\{0\})).$$

There is a diffeomorphism $f:Y_+\ra Y'$ such that

(1). We have that $f=C$ when restricted to $Y''\backslash{\rm int}({\rm im}(r''))=Y\backslash{\rm int}({\rm im}(r))$.

(2). We have that $f(r(\eta\times\{0\}))=r'(\eta'\times\{0\}$).

The diffeomorphism $f$ will induce a map
$$HM(f):HM(Y_+|r(R\times\{0\}))\ra HM(Y'|r(R\times\{0\})).$$

\brem
The first property actually implies that any two such $f$ would be isotopic to each other by lemma \ref{lem_unique_isotopy_class}. The second property ensures that $f$ will induce a map between monopole Floer homologies with local coefficients. One might have noticed that in the above construction, the result of cuting twice and re-gluing using $(\psi^{C})^{-1}$ and $\varphi^C\circ\psi^C$ is the same as cutting once and re-glue using $\varphi^C$. Yet the second property of $f$ above is the reason why we need not only $\varphi^C$ but also $\psi^C$.
\erem

The canonical map is defined as
\bdefn\label{defn_caninical_map_for_the_same_genus}
With the above notations, the canonical map
$$\Phi_{\mathcal{D},\mathcal{D}'}:SHM(\mathcal{D})\ra SHM(\mathcal{D}')$$
is defined as
$$\Phi_{\mathcal{D},\mathcal{D}'}=HM(f)\circ HM(X_+)\circ (HM(X_-))^{-1}.$$
\edefn

\section{Handle gluing maps and cancelations}
\subsection{Prelimilary discussions}
To start with, we first introduce the definition of contact handle attachments. The following definition is from Juh\'asz and Zemke \cite{juhasz1803contact}. (Also in Giroux \cite{giroux1991convexite} or Ozbagci \cite{ozbagci2011contact} or Baldwin and Sivek \cite{baldwin2016contact}.)

\bdefn\label{defn_contact_handle_attachment}
Suppose $(M,\ga)$ is a balanced sutured manifold. A {\it $3$-dimensional contact handle attachment of index $k$}, where $k\in\{0,1,2,3\}$, is a quadruple $h=(\phi,S,D^3,\delta)$. Here $D^3$ is a standard tight contact $3-$ball with $\delta$ being the dividing set on $\partial{D}$. Also $S\subset\partial{D}^3$ is a $2$-submanifold of $\partial{D}^3$ and 
$$\phi:S\ra \partial M$$
is the gluing diffeomorphism. The pair $(S,\phi)$ has different description for different index $k$:

(1). When $k=0$, $S=\emptyset$.

(2). When $k=1$, $S$ is the disjoint union of two disks. Each disk intersects the dividing set $\delta$ in a simple arc.

(3). When $k=2$, $S$ is an annulus on $\partial{D}^3$ and it intersects the dividing set $\delta$ in two simple arcs, and each simple arc represents a non-trivial class in $H_1(S,\partial{S})$.

(4). When $k=3$, $S=\partial{D}^3$.

\begin{figure}[h]
\centering
\begin{overpic}[width=3.0in]{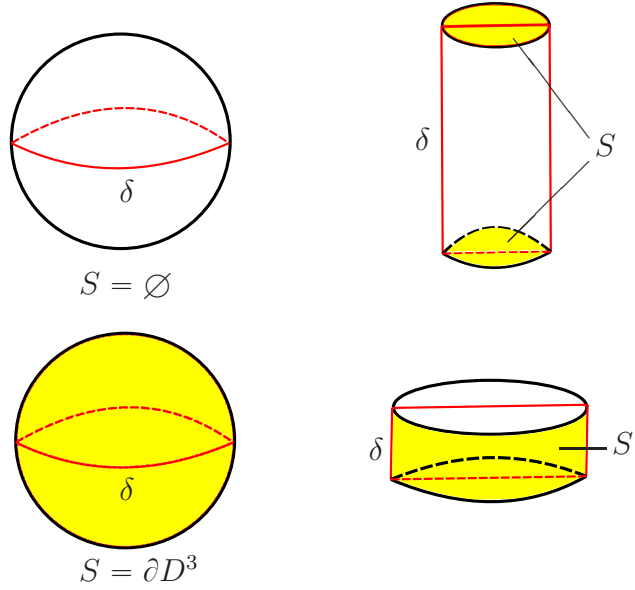}
	\put(12,45){$S=\emptyset$}
	\put(12,-5){$S=\partial{D}^3$}
	\put(19,61){$\delta$}
	\put(19,9){$\delta$}
	\put(70,70){$\delta$}
	\put(62,16){$\delta$}
	\put(100,72){\line(-3,4){13}}
	\put(100,68){\line(-1,-1){14}}
	\put(101,69){$S$}
	\put(103,18){\line(-1,0){8}}
	\put(104,17){$S$}
\end{overpic}
\caption{Contact handle attachment. Top left: $0$-handle. Top right: $1$-handle. Bottom right: $2$-handle. Bottom left: $3$-handle.}\label{fig_contact_handle_attachment}
\end{figure}

Furthermore, if we set $\delta_1=\delta\cap S$ and $\delta_2=\delta\backslash\delta_1$, then in any case we shall require that $\phi(\delta_1)\subset \ga\subset \partial{M}$. The new dividing set $\ga'$ for the new sutured manifold $M'=M\cup_{\phi}B$ is
$$\ga'=(\ga\backslash\phi(\delta_1))\cup(\delta_2).$$
See figure \ref{fig_contact_handle_attachment}.
\edefn

Before constructing handle gluing maps, let us first look at a special construction.

\begin{figure}[h]
\centering
\begin{overpic}[width=3.0in]{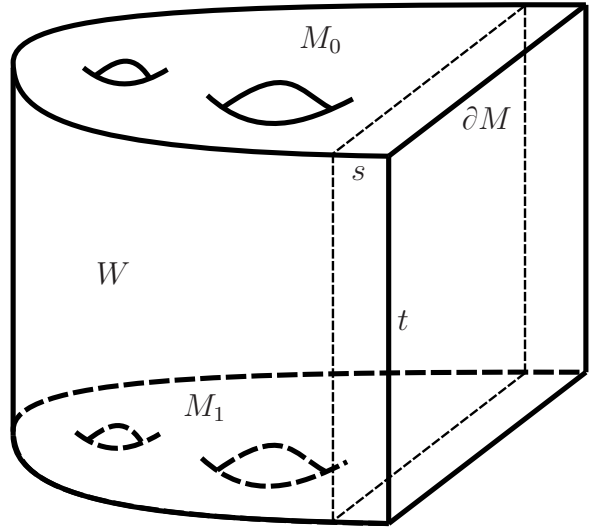}
\put(30,20){$M_1$}
\put(50,84){$M_0$}
\put(78,70){$\partial{M}$}
\put(59,61){$s$}
\put(15,43){$W$}
\put(67,35){$t$}
\end{overpic}
\vspace{0.05in}
\caption{The special type of cobordism $W$.}\label{fig_special_cobordism}
\end{figure}

\bprop\label{prop_map_from_special_cobordism}
Suppose $(M_0,\ga)$ and $(M_1,\ga)$ are two balanced sutured manifolds so that $\partial{M_0}=\partial{M_1}$ and the sutures are also identical. Suppose $W$ is a smooth compact oriented $4$-manifold so that $W$ can also be viewed as a manifold with conners: the boundary $\partial{W}$ consists of two horizontal parts $-M_0$ and $M_1$ as well as a vertical part $\partial{M}_0\times [0,1]$. The two parts $-M_0$ and $\partial{M}_0\times [0,1]$ meet in the conner $\partial{M}_0\times \{0\}$. The two parts $M_1$ and $\partial{M}_0\times [0,1]$ meet in the conner $\partial{M}_0\times \{1\}$. See figure \ref{fig_special_cobordism}. Then we can define a morphism between canonical modules:
$$F_{W}:\underline{\rm SHM}(M_0,\ga)\ra \underline{\rm SHM}(M_1,\ga).$$

\eprop
\bpf
Suppose $T$ is an auxiliary surface for $(M_0,\ga)$ and $f:\partial{T}\ra \ga$ is the map gluing $T$ to $(M_0,\ga)$. Let 
$$\widetilde{M}_0=M_0\mathop{\cup}_{f\times id}T\times[-1,1]$$
be the pre-closure and $\partial{\widetilde{M}}_0=R_+\cup R_-.$ Suppose $h:R_+\ra R_-$ is a diffeomorphism, we can use $h$ to glue $R_+\times [-1,1]$ to $\widetilde{M}_0$ and get a closure $Y_0$. Suppose $\eta$ is a non-separating curve on $R=R_+\times\{0\}$, we get a marked closure $\mathcal{D}_0=(Y_0,R,r,m_0,\eta)$ for $(M_0,\ga)$. Since the boundaries of $M_0$ and $M_1$ are identified, we can use the same auxiliary data $(T,f,h,\eta)$ to get a marked closure $\mathcal{D}_1=(Y_1,R,r,m_1,\eta)$.

There is a natural way to construct a cobordism from $Y_0$ to $Y_1$ out of $W$. Use $f\times id\times id$ to glue $T\times [-1,1]\times [0,1]$ to $A(\ga)\times[0,1]\subset\partial{M}_0\times [0,1]\subset \partial{W}$, and use $(id\cup h)\times id$ to glue $(R\times [0,1])\times [0,1]$ to the result of the first gluing. Finally we get a cobordism $\widehat{W}$ from $Y_0$ to $Y_1$. We have a map
$$HM(\widehat{W}): SHM(\mathcal{D}_0)\ra SHM(\mathcal{D}_1).$$

We claim that this map will induce a morphism between canonical modules. We only prove here that the cobordism map constructed above commutes with the canonical map $\Phi^g$, and the commutativity with $\Phi^{g,g+1}$ would follow with a similar argument. To proceed, suppose $\mathcal{D}_0'=(Y_0',R',r',m'_0,\eta')$ is another marked closure for $(M_0,\ga)$, obtained in a similar way as above, $\mathcal{D}_1'=(Y_1',R',r',m'_1,\eta')$ is the corresponding marked closure for $(M_1,\ga)$, and $\hat{W}'$ is the corresponding cobordism from $Y_0'$ to $Y_1'$. Then we need to show that the following diagram commutes up to multiplication by a unit:
\begin{equation*}
\xymatrix{
SHM(\mathcal{D}_0)\ar[rr]^{\Phi^g_{\mathcal{D}_0,\mathcal{D}_0'}}\ar[dd]^{HM(\widehat{W})}&&SHM(\mathcal{D}_0')\ar[dd]^{HM(\widehat{W}')}\\
&&\\
SHM(\mathcal{D}_1)\ar[rr]^{\Phi^g_{\mathcal{D}_1,\mathcal{D}_1'}}&&SHM(\mathcal{D}_1')\\
}
\end{equation*}

By definition \ref{defn_caninical_map_for_the_same_genus}, the canonical map $\Phi^g_{\mathcal{D}_0,\mathcal{D}_0'}$ is constructed as follows: the identity on $M_0$ can be extend to a diffeomorphism 
$$C_0: Y_0\backslash{\rm int}({\rm im}(r))\ra Y_0'\backslash{\rm int}({\rm im}(r')),$$
and there are maps $\varphi^{C_0}$ and $\psi^{C_0}$ and we can decompose them into composition of Dehn twists:

$$\varphi^C\circ\psi^C\sim D^{e_1}_{a_1}\circ...\circ D^{e_u}_{a_u},~\varphi^C\circ\psi^C\sim D^{e_{u+1}}_{a_{u+1}}\circ...\circ D^{e_v}_{a_v}.$$

For simplicity, we assume here that all $e_i=1$ (the general case follows from a similar argument) and the map $\Phi^g_{\mathcal{D}_0,\mathcal{D}_0'}$ is then induced by a cobordism $W_0$ obtained from $Y_1\times [0,1]$ attaching $4$-dimensional $2$-handles along curves $a_1,...,a_n\subset Y_1\times \{1\}.$ 

Since the two manifolds $M_0$ and $M_1$ have identical boundary: $\partial{M}=\partial{M}_1$, for constructing the canonical map $\Phi^g_{\mathcal{D}_1,\mathcal{D}_1'}$, we can chose a diffeomorphism

$$C_1: Y_1\backslash{\rm int}({\rm im}(r))\ra Y_1'\backslash{\rm int}({\rm im}(r'))$$

so that $C_1$ restrict to identity on $M_1$ and also
$C_0=C_1$ outside ${int}(M_0)$ and $int(M_1)$. Hence we have
$$\varphi^{C_1}=\varphi^{C_0},~\psi^{C_1}=\psi^{C_0}.$$
This means that the canonical map $\Phi^g_{\mathcal{D}_1,\mathcal{D}_1'}$ is induced by a cobordism $W_1$ which is obtained by attaching $4$-dimensional $2$-handles to $Y_1\times[0,1]$ along the same set of curves $a_1,...,a_n\subset Y_1\times\{1\}$.

Now the commutativity of the diagram is equivalent to
\begin{equation}\label{eq_5}
HM(W_1)\circ HM(\widehat{W})=HM(\widehat{W}')\circ HM(W_0).
\end{equation}
From the next lemma (lemma \ref{lem_handle_decomposition_of_cobordism}) we can view $\hat{W}$ as obtained from $M_0\times[0,1]$ by attaching $4$-dimensional handles $h^4_1,...,h^4_m$ to ${\rm int}(M_0)\times \{1\}$, and then $\hat{W}'$ is obtained from $Y_0\times[0,1]$ by attaching the same set of $4$-dimensional handles $h^4_1,...,h^4_m$ to ${\rm int}(M_0)\times \{1\}\subset Y_0'\times \{1\}$. So to prove the equality (\ref{eq_5}), it is enough to prove that the set of handles  $h^4_1,...,h^4_m$ and the set of $2$-handles attached along $a_1,...,a_n\subset Y_0\times\{1\}$, which are coming from the construction of canonical maps between closures, can commute with each other. But this is obvious: $h^4_1,...,h^4_m$ are attached to ${\rm int}(M_0)\times \{1\}\subset Y_0\times \{1\}$ and the curves $a_1,...,a_n$ are inside ${\rm int}({\rm im}(r))\times \{1\}\subset Y_0\times \{1\}$ and
$${\rm int}(M_0)\cap {\rm int}({\rm im}(r))=\emptyset.$$
\epf

\blem\label{lem_handle_decomposition_of_cobordism}
Suppose $(M_0,\ga)$, $(M_1,\ga)$ and $W$ are defined as in proposition \ref{prop_map_from_special_cobordism}. Then $W$ is diffeomorphic to a $4$-manifold obtained from $M_0\times [0,1]$ by attaching some $4$-dimensional handles to ${\rm int}(M_0)\times \{1\}$. 
\elem

\bpf
We can assume a neighborhood $N$ of the vertical boundary part $\partial{M_1}\times [0,1]$ of $W$ is identified with $\partial{M}_0\times [-1,0]_s\times [0,1]_t$ so that the vertical boundary part is $\partial{M}_0\times\{0\}\times[0,1]$. We can choose a smooth function $f:W\ra [0,1]$ so that

$$f(-M_0)=0,~f(M_1)=1,~f(\partial{M_0}\times[-1,0]\times\{t\})=t.$$
Perturb $f$ a little bit so that $f$ is Morse and there is no critical points of $f$ near $\partial{W}\subset W$. Such perturbation exists since the set of Morse functions is dense in the space of smooth functions and $f$ has already been Morse near the boundary $\partial{W}\subset W$, and having no critical points there. Then $f$ induces the desired handle decomposition.
\epf
\brem
In \cite{juhasz2012naturality} Juh\'asz and Thurston also proved that $0$- and $4$- handle attachments can be avoided.
\erem

Suppose we are giving a smooth compact oriented $4$-manifold with boundary $W$, and let $S\subset \partial{W}$ be a closed oriented surface surface which separates $\partial{W}$ into two parts. Let $M_1$ and $M_2$ be the closures of those two parts with orientations so that 
$$\partial{M_1}=\partial{M}_2=S, ~-M_1\cup M_2=\partial{W}.$$
Suppose $\ga\subset S$ is a collection of oriented simple closed curves so that $(M_1,\ga)$ and $(M_2,\ga)$ are all balanced sutured manifolds. We can view $W$ as a cobordism from $(M_1,\ga)$ to $(M_2,\ga)$. An adaption of lemma \ref{lem_handle_decomposition_of_cobordism} shows that $W$ is actually diffeomorphic to the $4$-manifolds obtained from $M_1\times[0,1]$ by attaching some $4$-dimensional handles along ${\rm int}(M_1)\times\{1\}.$ Hence just as in the proof of proposition \ref{prop_map_from_special_cobordism} we can also have a map between sutured monopole Floer homologies of $(M_1,\ga)$ and $(M_2,\ga)$. Sometimes it is more convenient to use this setting so we will give this a name:

\bdefn
Under the above settings, we call $W$ a cobordism with {\it sutured surface} $(S,\delta)$, from $(M_1,\ga_1)$ to $(M_2,\ga_2)$. The collection of oriented simple curves $\delta$ on $S$ is called a {\it suture}.
\edefn

\subsection{Constructions of handle gluing maps}
The definition of handle attaching maps for attaching $0$- and $1$-handles in \cite{baldwin2016contact} are kind of straightforward and are summarized as follows.

\bdefn\label{defn_zero_and_one_handle_gluing_map}
Suppose $(M,\ga)$ is a balanced sutured manifold and $h$ is a $0$ or $1$-handle attached to $(M,\ga)$ and results in a new balanced sutured manifold $(M',\ga')$. We can use auxiliary data $T$ and $f$ to produce $\widetilde{M}$ which is a pre-closure of $(M,\ga)$. Then $\widetilde{M}$ is also a pre-closure of $(M',\ga')$. Hence we can use the same auxiliary data $R,h,\eta$ to get marked closures $\mathcal{D}'=(Y,R,r,m',\eta)$ and $\mathcal{D}=(Y,R,r,m=m'|_{M},\eta)$ for $(M',\ga')$ and $(M,\ga)$ respectively. The {\it handle attaching map} 
$$C_h:\underline{\rm SHM}(-M,-\ga)\ra\underline{\rm SHM}(-M',-\ga')$$
for the contact handle $h$ is then defined to be the map between projective transitive systems induced by the product cobordism $(-Y)\times[0,1]$ (or simply just the identity map).
\edefn

\begin{figure}[h]
\centering
\begin{overpic}[width=5.0in]{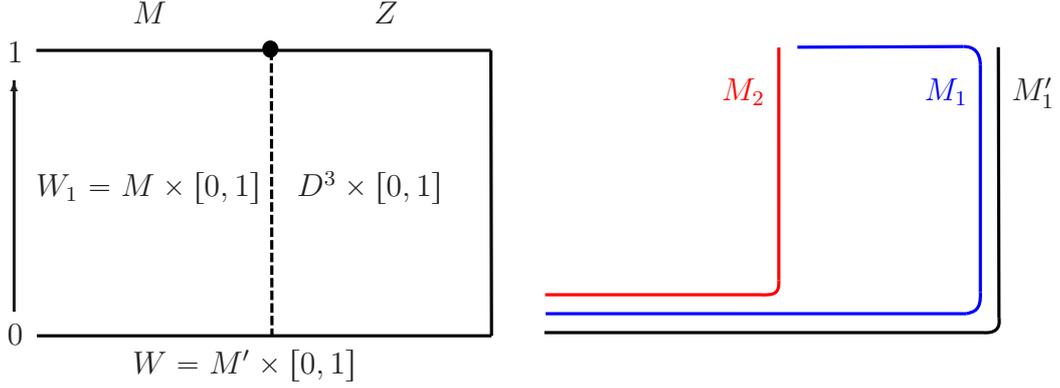}
	\put(-3,-0.5){$0$}
	\put(-2.3,3){\vector(0,1){24}}
	\put(-3,29){$1$}
	\put(10,33){$M$}
	\put(35,33){$Z$}
	\put(10,-3.5){$W=M'\times[0,1]$}
	\put(0,15){$W_1=M\times[0,1]$}
	\put(27,15){$D^3\times[0,1]$}
	\put(71,25){\color{red}$M_2$}
	\put(92,25){\color{blue}$M_1$}
	\put(101,25){$M_1'$}
\end{overpic}
\vspace{0.05in}
\caption{The product $W,W_1$, and the sutured manifolds $M_1,M_2,M_1'$.}\label{fig_cobordism}
\end{figure}

The paper \cite{baldwin2016contact} also discussed on the handle attaching maps for contact $2$- and $3$-handles. But we want to introduce somewhat different definitions for our convenience. Suppose we are attaching a $2$- or $3$-handle $h=(\phi,S,D^3,\delta)$ and result in $(M',\ga')$. Let $Z=M'\backslash{\rm int}(M)$. The idea is that when we turn $h$ up-side down, we will get a $1$- or $0$- handle as a result. To turn $h$ up-side down, we shall consider the manifold $W=M'\times[0,1]$. We can chose the surface
$S=\partial{M}\times\{0\}\subset M_1\times\{0\}\subset\partial{W}$ and the suture $\ga=\ga\times\{0\}\subset S$. Then $W$ can be viewed as a cobordism with sutured surface $(S,\ga)$, from $(M,\ga)$ to
$$(M_1,\ga)=(M'\cup\partial{M'}\times[0,1]\cup Z,\ga).$$
In this case, $Z$ is attached to $\partial{M}'\times\{1\}$ and can be viewed as a $0$- or $1$-handle $h^{\vee}$. Let $M_1'=M'\cup\partial{M}'\times[0,1]$ and $\ga_1'=\ga'\times\{1\}\subset \partial{M}_1'$. See figure \ref{fig_cobordism}.

\vspace{0.1in}

Now we have a handle attaching map
$$C_{h^{\vee}}:\underline{SHM}(-M_1',-\ga'_1)\ra\underline{SHM}(-M_1,-\ga)$$
which is an isomorphism by definition \ref{defn_zero_and_one_handle_gluing_map}. Proposition \ref{prop_map_from_special_cobordism} induces a map
$$F_{-W}:\underline{SHM}(-M,-\ga)\ra\underline{SHM}(-M_1,-\ga),$$
so we only need a map
$$\Psi: \underline{SHM}(-M_1',-\ga'_1)\ra\underline{SHM}(-M',-\ga').$$
This map seems to be obvious, since 
$$(M'_1,\ga'_1)=(M'\cup \partial{M}'\times[0,1],\ga'\times\{1\})$$
is just gotten from $(M',\ga')$ by attaching a collar of the boundary. This can be made precise as follows. When closing up $(M_1',\ga_1')$, we choose an auxiliary surface $T_1'$ and glue $T_1'\times[-1,1]$ to $(M_1',\ga_1')$ along $A(\ga'_1)=\ga'_1\times[-1,1]$ by a map
$$f:\partial{T_1'}\ra \ga'_1.$$
Let $\widetilde{M}'_1=M_1'\cup T_1'\times[-1,1]$ and suppose
$$\partial{\widetilde{M}_1'}=R_+'\cup R_-',$$
so that $R_{\pm}'$ contains $R_{\pm}(\ga')$. If we choose an orientation preserving diffeomorphism $h:R_+'\ra R_-'$, then we can glue $R_+'\times[0,1]$ to $Y'_1$ and get a marked closure $\mathcal{D}'_1=(Y'_1,R'_+,r_1,m_1,\delta_1)$ for $(M_1',\ga_1')$. 

\begin{figure}[h]
\centering
\begin{overpic}[width=4.0in]{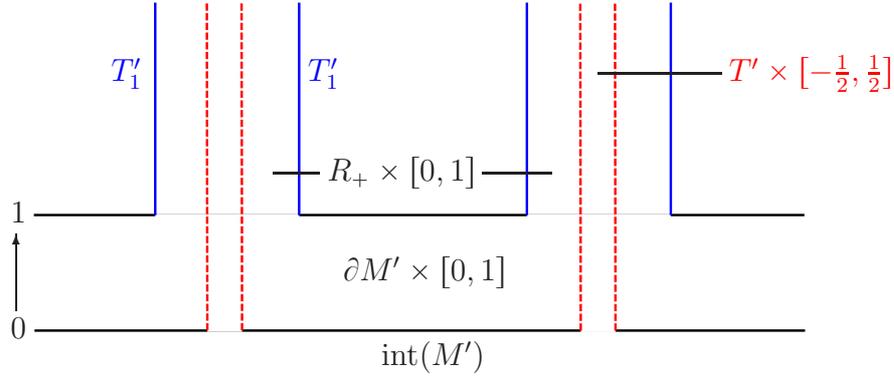}
	\put(-3,-0.5){$0$}
	\put(-2.3,3){\vector(0,1){10}}
	\put(-3,14.5){$1$}
	\put(10,33){\color{blue}$T_1'$}
	\put(35.5,33){\color{blue}$T_1'$}
	\put(90,33){\color{red}$T'\times[-\frac{1}{2},\frac{1}{2}]$}
	\put(89,34){\line(-1,0){16}}
	\put(40,7){$\partial M'\times[0,1]$}
	\put(38,20){$R_+\times[0,1]$}
	\put(37,21){\line(-1,0){6}}
	\put(58,21){\line(1,0){9}}
	\put(45,-4){${\rm int}(M')$}
\end{overpic}
\vspace{0.1in}
\caption{The collar $\partial{M}'\times[0,1]$ and the auxiliary surfaces $T_1',T'$.}\label{fig_collar_of_boundary}
\end{figure}

We want to show next that there is a canonical way to view $\mathcal{D}_1'$ as a closure of $(M',\ga')$. We can view the original collar $A(\ga')$ of $\ga'$ as identified with $\ga'\times[-1,1]$. We can get a new product neighborhood $A'(\ga')=\ga'\times[-\frac{1}{2},\frac{1}{2}]$. Now let 
$$T'=T_1'\cup_{f}\ga'\times[0,1],$$
where $f:\partial{T_1'}\ra \ga'=\ga_1'$ is the map defined as above and the annuli glued to $T_1'$ via $f$ are chosen to be $\ga'\times[0,1]\subset \partial{M'}\times[0,1]$. Then we can view $T'\times[-\frac{1}{2},\frac{1}{2}]$ as attached to $(M',\ga')$ along $A'(\ga')$. Let $\widetilde{M}'=M'\cup T'\times[-\frac{1}{2},\frac{1}{2}]$, we want to show that there is a canonical way to identify
\begin{equation}\label{eq_unique_isotopy_of_collar}
	\widetilde{M}'_1\backslash{\rm int}(\widetilde{M}')\cong (R_+'\sqcup R_-')\times[-1,1].
\end{equation}
See figure \ref{fig_collar_of_boundary}.

\vspace{0.1in}

Suppose $\partial\tilde{M}'=R_+\cup R_-$, then there is a canonical way up to isotopy to identify $R_{\pm}$ with $R_{\pm}'$. The $R_{\pm}(\ga')$ is identified with $R_{\pm}(\ga'_1)$ and $T_1'$ part is identified with itself. The rest of $R_{\pm}$ are product annuli of $\ga'$. Hence there is a canonical way to identify
\begin{equation}\label{eq_canonical_identification_on_boundary}
	\partial{(\widetilde{M}'_1\backslash{\rm int}(\widetilde{M}'))}\cong (R_+'\sqcup R_-')\times\{-1,1\}.
\end{equation}
Note it is obvious that 
$$(\widetilde{M}'_1\backslash{\rm int}(\widetilde{M}'))\cong (R_+'\sqcup R_-')\times[0,1],$$
so by lemma \ref{lem_unique_isotopy_class}, there is a unique isotopy class of diffeomorphisms which restrict to the canonical identification (\ref{eq_canonical_identification_on_boundary}) on the boundary. Hence there is a well defined map
$$\Psi: \underline{\rm SHM}(-M_1',-\ga'_1)\ra\underline{\rm SHM}(-M',-\ga')$$
induced by the identity on $HM(-Y_1'|-R_1';\Gamma_{-\eta_1'})$.

\bdefn\label{defn_two_and_three_handle_gluing_map}
Suppose $h$ is a $2$- or $3$-handle attached to $(M,\ga)$ and $(M',\ga')$ is the resulting balanced sutured manifold. Suppose $C_{h^{\vee}}$, $F_{-W}$ and $\Psi$ are defined as above, then we define the contact handle attaching map as
$$C_h=\Psi^{-1}\circ C_{h^{\vee}}^{-1}\circ F_{W}:\underline{SHM}(-M,-\ga)\ra \underline{SHM}(-M',-\ga').$$
\edefn

Actually the handle gluing maps constructed above are the same as what are done in Baldwin and Sivek's paper \cite{baldwin2016contact}.

\bprop\label{prop_equivalence_of_the_definition}
The gluing maps in definition \ref{defn_two_and_three_handle_gluing_map} are equivalent to those constructed by Baldwin and Sivek in \cite{baldwin2016contact}. 
\eprop
\bpf

For $2$- or $3$-handles, suppose $(M,\ga)$ is the original sutured manifold and $(M',\ga')$ is the result of attaching a contact handle $h=(\phi,S,D^3,\delta)$. Suppose $(M_1,\ga)$ and $(M'_1,\ga_1')$ are constructed as in definition \ref{defn_two_and_three_handle_gluing_map}.

Suppose $W=M'\times[0,1]$ is the product, we can view it as a cobordism with sutured surface $(S=\partial{M}\times\{0\},\ga)$. When doing closing up along $(S,\ga)$, we get two marked closures $\mathcal{D}=(Y,R,r,m,\eta)$ and $\mathcal{D}_1=(Y_1,R,r,m_1,\eta)$ for $(M,\ga)$ and $(M_1,\ga)$ respectively and a cobordism $\widehat{W}$ from $Y$ to $Y_1$ which induces the map $F_{-W}$.  From definition \ref{defn_zero_and_one_handle_gluing_map} and the construction of $\Psi$ in definition \ref{defn_two_and_three_handle_gluing_map}, we can see that $\mathcal{D}_1$ is also a marked closure of $(M',\ga')$. Thus handle attaching map is just induced by the cobordism $\widehat{W}$.

Let $W_1=M\times[0,1]\subset W$ and we can view $W_1$ as a special cobordism with sutured surface $(S_1,\ga)$. Let $M_2=M\times\{1\}\cup\partial{M}\times[0,1]$. See figure \ref{fig_cobordism}. By doing a suitable closing up along $(S_1,\ga)$, we can get two marked closures $\mathcal{D}$ (the same as above) and $\mathcal{D}_2=(Y_2,R,r,m_2,\eta)$ for $(M,\ga)$ and $(M_2,\ga)$ respectively and a cobordism $\widehat{W}_1$ from $Y$ to $Y_2$. 

Recall $h=(\phi,S,D^3,\delta)$ is the contact handle attached to $M$. Then $W$ can be viewed as obtained from $W_1$ by attaching $D^3\times[0,1]$ through the map
$$\phi\times id:S\times[0,1]\ra \partial{M}\times[0,1]\subset M_2\subset \partial{W}_1.$$
Accordingly, $\widehat{W}$ can be viewed as obtained from $\widehat{W}_1$ by attaching $D^3\times[0,1]$ through the map
$$\phi\times id:S\times[0,1]\ra \partial{M}\times[0,1]\subset M_2\subset Y_2\subset\partial{\widehat{W}}_1.$$
This $D^3\times[0,1]$ now becomes a $4$-dimensional handle. Let $h^4=(\phi\times id, S\times[0,1],D^3\times[0,1])$ be this $4$-dimensional handle.

If $h$ is a $3$-dimensional $2$-handle, then $h^4$ is a $4$-dimensional $2$-handle. Recall for a $2$-handle $h$, $S$ is an annulus. Suppose $\al\subset {\rm int}(S)$ is the core of $S$, i.e., 
$$[\al]=\pm1\subset H_1(S)\cong\intg,$$
then we can view $h$ as a $2$-handle attached along the curve $\phi(\al)\subset \partial{M}$, and view $h^4$ as attached along the curve 
$$\phi(\al)\times\{\frac{1}{2}\}\subset Y_2\subset  \partial{\widehat{W}}_1.$$
Note that $\widehat{W}_1$ is actually diffeomorphic to $Y\times[0,1]$ so such a $4$-dimensional $2$-handle attachment actually corresponds to a Dehn surgery on $Y$ along the curve $\phi(\al)\subset \partial{M}\subset Y$. The slope of the Dehn surgery can be compute from the framing of the $4$-dimensional $2$-handle attached and it is a $0$ Dehn surgery with respect to $\partial{M}\times\{\frac{1}{2}\}$ surface framing. If $\partial{M}\times[0,1]$ is equipped with $I$-invariant contact surface so that $\partial{M}$ is convex and $\ga$ is the dividing set, then since $\phi(\al)$ intersects $\ga$ twice we know that when we realize $\phi(\al)$ as a Legendrian curve intersecting $\ga$ twice, then the $0$ surface framing corresponds to the contact $+1$ framing. This is exactly the case in \cite{baldwin2016contact}. Hence the two gluing maps are the same.

If $h$ is a $3$-handle, then $h^4$ is also a $3$-handle. Now $\hat{W}$ is diffeomorphic to the result of attaching a $4$-dimensional $3$-handle to $Y\times\{1\}\subset Y\times[0,1]$ hence it is also diffeomorphic to the result of gluing a $4$-dimensional $1$-handle to $Y_2\times[0,1]$ along two points in $Y_2\times\{0\}$. This is exactly the same as in \cite{baldwin2016contact} so the two gluing maps are the same.
\epf

Above discussions will also lead to an equivalent definition for the handle gluing maps for $2$- and $3$-handles, which would be more useful in later discussions. 

\bdefn\label{defn_two_handle_gluing_map2}
Suppose $(M,\ga)$ is a balanced sutured manifold and $h$ is a $2$-handle attached to $(M,\ga)$ along the curve $\al\subset \partial{M}$ and results in a new balanced sutured manifold $(M',\ga')$. Suppose we use auxiliary data $T,f$ to get a pre-closure $\widetilde{M}$ for $(M,\ga)$, then we can do a $0$-Dehn surgery along a curve $\be\subset{\rm int}(M)$, which is isotopic to $\al$, with respect to $\partial{M}$-surface framing to get a pre-closure $\widetilde{M}'$ for $(M',\ga')$. Since the surgery is supported in a neighborhood of $\be$, we have $\partial{\widetilde{M}}=\partial{\widetilde{M}}'$ and hence we can use the same auxiliary data $R,h,\eta$ to get marked closures $\mathcal{D}'=(Y,R,r,m,\eta)$ and $\mathcal{D}=(Y',R,r,m',\eta)$ for $(M,\ga)$ and $(M',\ga')$ respectively. We can form a cobordism $\widehat{W}$ from $Y$ to $Y'$ obtained by attaching a $0$-framed $4$-dimensional $2$-handle to $Y\times[0,1]$ along $\be\subset Y\times\{1\}$. The {\it handle gluing map} 
$$C_h:\underline{\rm SHM}(-M,-\ga)\ra\underline{\rm SHM}(-M',-\ga')$$
for a $2$-handle $h$ is defined to be the map induced by the cobordism $-\widehat{W}$.
\edefn

\bdefn\label{defn_three_handle_gluing_map2}
Suppose $(M,\ga)$ is a balanced sutured manifold and $h$ is a $2$-handle attached to $(M,\ga)$ along a sphere $S\subset \partial{M}$ and results in a new balanced sutured manifold $(M',\ga')$. Suppose we use auxiliary data $T,f$ to get a pre-closure $\widetilde{M}$ for $(M,\ga)$, then we can do a cut and paste surgery along a sphere $S'\subset{\rm int}(M)$, which is isotopic to $S$, to get a pre-closure $\widetilde{M}'$ for $(M',\ga')$. Since the surgery is supported in a neighborhood of $S'$, we have $\partial{\widetilde{M}}=\partial{\widetilde{M}}'$ and hence we can use the same auxiliary data $R,h,\eta$ to get marked closures $\mathcal{D}'=(Y,R,r,m,\eta)$ and $\mathcal{D}=(Y',R,r,m',\eta)$ for $(M,\ga)$ and $(M',\ga')$ respectively. We can form a cobordism $\widehat{W}$ from $Y$ to $Y'$ obtained by attaching a $4$-dimensional $3$-handle to $Y\times[0,1]$ along $S'\subset Y\times\{1\}$ corresponding to the cut and paste surgery. The {\it handle gluing map} 
$$C_h:\underline{\rm SHM}(-M,-\ga)\ra\underline{\rm SHM}(-M',-\ga')$$
for a $3$-handle $h$ is defined to be the map induced by the cobordism $-\widehat{W}$.
\edefn

\subsection{Basic properties of handle attaching maps}
For a special pair of handles, we can cancel them both topologically and for cobordism maps.

\blem\label{lem_zero_one_cancelation}
Suppose $(M,\ga)$ is a balanced sutured manifold, $h=(\phi,S,D^3,\delta)$ is a $0$-handle and $h'=(\phi',S',D^{3}{}',\delta')$ is a $1$-handle such that the attaching map $\phi'$ maps one component of $S'$ to $\partial{M}$ and the other component to $\partial{D}^3$. Let $(M',\ga')$ be the resulting balanced sutured manifold, then there is a diffeomorphism $\psi_{01}:(M,\ga)\ra (M',\ga')$ so that

(1). The map $\psi_{01}$ restricts to identity outside a neighborhood of $\phi'(S')\cap \partial{M}$.

(2). The map $\psi_{01}$ is isotopic to the inclusion map $M\hookrightarrow M'$

Further more, we have
$$C_{h'}\circ C_{h}=\underline{\rm SHM}(\psi_{01}):\underline{\rm SHM}(-M,-\ga)\ra \underline{\rm SHM}(-M',-\ga').$$
\elem

\bpf
The two handles form a pair of handles which can be canceled topologically, so we can easily find such a diffeomorphism $\psi_{01}:(M,\ga)\ra(M',\ga')$ satisfying the two conditions above.

From definition \ref{defn_zero_and_one_handle_gluing_map}, we know that a marked closure $\mathcal{D}'=(Y,R,r,m',\eta)$ will induce a marked closure $\mathcal{D}=(Y,R,r,m=m'|_{M},\eta)$ and the composition $C_{h'}\circ C_{h}$ is induced by the identity map 
$$id:SHM(-\mathcal{D})\ra SHM(-\mathcal{D}').$$

Let us now describe the map $\underline{\rm SHM}(\psi_{01})$. If we fix the same closure $\mathcal{D}'=(Y,R,r,m',\eta)$, then we can get a closure 
$$\widetilde{\mathcal{D}}=(Y,R,r,m'\circ \psi_{01},\eta)$$
for $(M,\ga)$ and the map $\underline{SHM}(\psi_{01})$ is the induced by the map
$$id:SHM(-\widetilde{\mathcal{D}})\ra SHM(-\mathcal{D}').$$
To prove the lemma, we need to show the commutativity up to multiplication by a unit of the following diagram, by lemma \ref{lem_compare_canonical_maps}:
\begin{equation*}
\xymatrix{
SHM(\mathcal{-D})\ar[rrdd]^{id}\ar[rr]^{\Phi^g_{-\mathcal{D},-\widetilde{\mathcal{D}}}}&&SHM(-\widetilde{\mathcal{D}})\ar[dd]^{id}\\
&&\\
&&SHM(-\mathcal{D}')\\
}
\end{equation*}

Now since $m$ and $m'\circ \psi_{01}$ are isotopic in $Y$ so that the isotopy is identity outside a neighborhood of $\phi''(S'')\subset{\rm int}(Y\backslash{\rm int}({\rm im}(r)))$ and the two marked closures have the same $r$ map, the canonical map from $SHM(-\mathcal{D})$ to $SHM(-\widetilde{\mathcal{D}})$ is just the identity map by definition \ref{defn_caninical_map_for_the_same_genus}, so the above diagram indeed commute, and we must have
$$C_{h'}\circ C_{h}=\underline{\rm SHM}(\psi_{01}):\underline{\rm SHM}(-M,-\ga)\ra \underline{\rm SHM}(-M',-\ga').$$
\epf

\blem\label{lem_one_two_handle_cancelation}
Suppose $(M,\ga)$ is a balanced sutured manifold, $h=(\phi,S,D^3,\delta)$ is a $1$-handle attached to $(M,\ga)$ and $h'=(\phi',S',D^{3}{}',\delta')$ is a $2$-handle attached to $(M,\ga)\cup h$. Suppose $\al'\subset S'$ is the core of the annulus $S'$, i.e., the simple closed curve which generates $H_1(S)$. Suppose the attaching map $\phi'$ maps the core 
$\al'$ to a curve which is the union of two arcs
$$\phi'(\al')=c_1'\cup c_2'\subset \partial{(M\mathop{\cup}_{\phi}D^3)}.$$
Here $c_1'$ is an arc on $\partial{M}$ disjoint with the suture $\ga$ and $c_2'$ intersects the dividing set $\delta$ on $\partial{D}^3\backslash S$ twice, and  we shall require that $c_2'$ intersects each component (there are two components) of $\delta\backslash S$ once.

Suppose $(M_2,\ga_2)$ is the resulting manifold of attaching $h$ and $h'$, then there is a canonical isotopic class of diffeomorphisms $\psi_{12}: (M,\ga)\ra (M_2,\ga_2)$ so that 

(1). The map $\psi_{12}$ restricts to identity outside a neighborhood of $(\phi(S)\cup\phi'(S'))\cap\partial{M}\subset M$.

(2). The map $\psi_{12}$ is isotopic to the inclusion $M\hookrightarrow M'$.

Furthermore, we have
$$C_{h'}\circ C_{h}=\underline{\rm SHM}(\psi_{12}):\underline{\rm SHM}(-M,-\ga)\ra\underline{\rm SHM}(-M_2,-\ga_2).$$
\elem

\begin{figure}[h]
\centering
\begin{overpic}[width=4.0in]{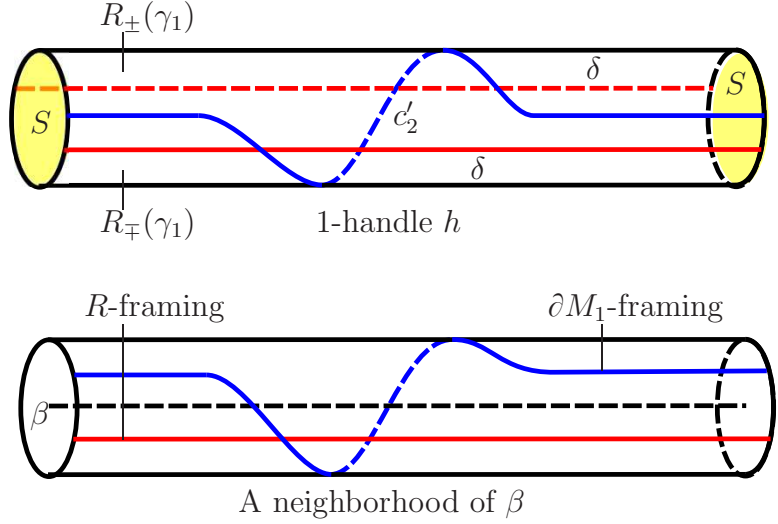}
	\put(40,33){$1$-handle $h$}
	\put(60,40){$\delta$}
	\put(15,54){\line(0,1){5}}
	\put(12,60){$R_{\pm}(\ga_1)$}
	\put(15,41){\line(0,-1){5}}
	\put(12,33){$R_{\mp}(\ga_1)$}
	\put(75,53){$\delta$}
	\put(50,47){$c_2'$}
	\put(3,46){$S$}
	\put(93,51){$S$}
	\put(30,-3.5){A neighborhood of $\be$}
	\put(3,8){$\be$}
	\put(15,6){\line(0,1){15}}
	\put(10,22){$R$-framing}
	\put(70,22){$\partial{M_1}$-framing}
	\put(77,15){\line(0,1){6}}
\end{overpic}
\vspace{0.05in}
\caption{Top: the $1$-handle $h$ and the part $c_2'$ of the core of the $2$-handle $h'$. Bottom: in a neighborhood of the curve $\be$, the longitudes of the two surface framings}\label{fig_one_two_cancelation}
\end{figure}

\bpf
The two handles can be canceled topologically so the map $\psi_{12}$ is easy to find.

Now suppose $(M_1,\ga_1)$ is the result of attaching the $1$-handle $h$. Then from definition \ref{defn_zero_and_one_handle_gluing_map} we know that a marked closure $\mathcal{D}_1=(Y,R,r,m_1,\eta)$ of $(M_1,\ga_1)$ will induce a marked closure $\mathcal{D}=(Y,R,r,m=m_1|_{M},\eta)$ of $(M,\ga)$ and the map $C_{h}$ is induced by the identity map 
$$id:SHM(-\mathcal{D})\ra SHM(-\mathcal{D}_1).$$

From definition \ref{defn_two_handle_gluing_map2}, there is a curve $\be\subset Y$ isotopic to $m_1(\phi(\al))\subset m_1(\partial{M}_1)\subset Y$ so that if we do a $0$-Dehn surgery with respect to $m_1(\partial{M}_1)$ surface framing, then the resulting manifold $Y_2$ is a closure of $(M_2,\ga_2)$. Now $\be\subset {\rm int}(m_1(M_1))\subset Y$ and the Dehn surgery can be supported in an arbitrarily small tubular neighborhood of $\be$. Hence the data for $r,R,\eta$ in $\mathcal{D}_1$ is not influenced by the Dehn surgery along $\be$ and we get a marked closure $\mathcal{D}_2=(Y_2,R,r,m_2,\eta)$ for $(M_2,\ga_2)$. As in definition \ref{defn_two_handle_gluing_map2}, the Dehn surgery corresponds to a cobordism $\widehat{W}$ from $Y_1$ to $Y_2$, obtained from $Y_1\times[0,1]$ by attaching a $0$-framed $4$-dimensional $2$-handle along $\be\times\{1\}\subset Y_1\times\{1\}$. So we have a cobordism
$$HM(-\widehat{W}):SHM(-\mathcal{D}_1)\ra SHM(-\mathcal{D}_2)$$
and this map induces $C_{h'}$.

Let us now describe the map $\underline{\rm SHM}(\psi_{12})$. If we fix the same closure $\mathcal{D}_2=(Y_2,R,r,m_2,\eta)$, then we can get a closure 
$$\widetilde{\mathcal{D}}=(Y_2,R,r,m_2\circ \psi_{12},\eta)$$
for $(M,\ga)$ and the map $\underline{\rm SHM}(\psi_{12})$ is the induced by the map
$$id:SHM(-\widetilde{\mathcal{D}})\ra SHM(-\mathcal{D}_2).$$
by lemma \ref{lem_compare_canonical_maps}, to finish the proof, we need to show the commutativity, up to multiplication by a unit, of the following diagram:
\begin{equation*}
\xymatrix{
SHM(\mathcal{-D})\ar[rrdd]^{HM(-\widehat{W})\circ id}\ar[rr]^{\Phi^g_{-\mathcal{D},-\widetilde{\mathcal{D}}}}&&SHM(-\widetilde{\mathcal{D}})\ar[dd]^{id}\\
&&\\
&&SHM(-\mathcal{D}_2)\\
}
\end{equation*}

Now let us describe $\Phi^g_{-\mathcal{D},-\widetilde{\mathcal{D}}}$ in details. The key observation is that under the condition of the lemma, we can isotope $\be$ into $R_+(\ga_1)$ or $R_-(\ga_1)$ and then into a curve $\be'\subset r_1(R\times \{t\})$ for any $t\in(-1,1)$. The reason is that from the hypothesis of the lemma we know that $\phi'(\al')=c_1'\cup c_2'$, $c_1'$ has already been contained in $R_+(\ga_1)$ or $R_-(\ga_1)$, and $c_2'$ can be isotoped into the same component within the $1$-handle $h$.

The surgery on $\be$ is $+1$ with respect to the contact framing and $0$ with respect to the $\partial{M_1}$-surface framing as discussed in the proof of proposition \ref{prop_equivalence_of_the_definition}. It is straightforward to see that after the isotopy, the surgery becomes a $\pm1$-surgery along $\be'$ with respect to the surface $r_1(R\times\{t\})$. See figure \ref{fig_one_two_cancelation}. If we go through contact framing again, since now $\be'$ does not intersect the dividing set on $r_1(R\times\{t\})$, we can see that it is a $+1$-surgery. When reverse the orientation to deal with $-\mathcal{D}_1$ and $-\widetilde{\mathcal{D}}_1$, it becomes a $-1$-surgery and hence corresponds to a positive Dehn twist. Hence from the definition \ref{defn_caninical_map_for_the_same_genus} for the canonical map, we know that the canonical map is induced by a cobordism $\hat{W}'$ which is obtained from $Y\times[0,1]$ by attaching a $4$-dimensional $2$-handle to along $\be'\subset Y\times\{1\}$. Yet $\widehat{W}$ and $\widehat{W}'$ are diffeomorphic since $\be$ and $\be'$ are isotopic and the framing of the handle gluing are also the same. Hence the above diagram indeed commute and we are done.
\epf

\begin{figure}[h]
\centering
\begin{overpic}[width=4.0in]{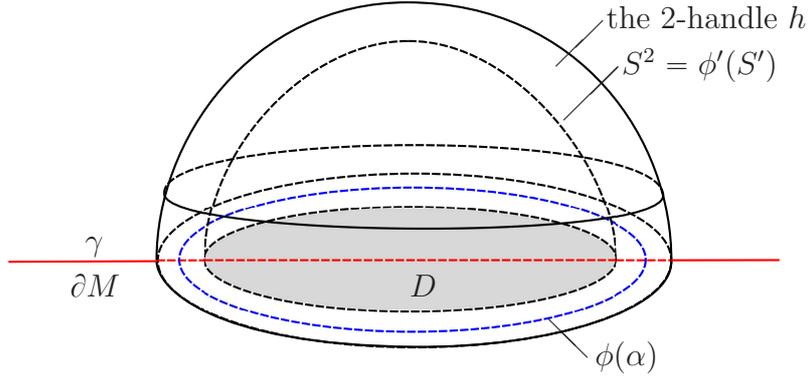}
	\put(10,13){$\ga$}
	\put(8,7){$\partial{M}$}
	\put(70,4){\line(1,-1){5}}
	\put(76,-2){$\phi(\al)$}
	\put(52,7){$D$}
	\put(72,30){\line(1,1){7}}
	\put(79.5,36){$S^2=\phi'(S')$}
	\put(71,37){\line(1,1){6}}
	\put(77.5,42){the $2$-handle $h$}
\end{overpic}
\vspace{0.05in}
\caption{The $2$-handle $h$ and the sphere $S^2$ along which $h'$ is attached.}\label{fig_two_three_cancelation}
\end{figure}

\blem\label{lem_two_three_handle_cancelation}
Suppose $(M,\ga)$ is a balanced sutured manifold, $h=(\phi,S,D^3,\delta)$ is a $2$-handle and $h'=(\phi',S',D^{3}{}',\delta')$ is a $3$-handle both attached to $(M,\ga)$. If $\al\subset {\rm int}(S)$ is a curve which represents a generator of $H_1(S)$, then we shall require that $\al$ is mapped to a curve on $\partial{M}$ which intersects $\ga$ twice and bounds a disk $D$ on $\partial{M}$. Hence a retraction of this disk $D$ union with one component of $\partial{D}^3\backslash S$ will become a new sphere boundary $S^2$ of the resulting manifold $(M_1,\ga_1)$ of gluing $h$ to $(M,\ga)$. We shall require that the attaching map $\phi'$ maps $S'=\partial{D}^3{}'$ to $S^2$. See figure \ref{fig_two_three_cancelation}.

Suppose $(M_2,\ga_2)$ is the resulting manifold of attaching $h$ and $h'$, then there is a canonical isotopic class of diffeomorphisms $\psi_{23}: (M,\ga)\ra (M_2,\ga_2)$ so that

(1). The map $\psi_{23}$ restricts to identity outside a neighborhood of $D\subset M$.

(2). The map $\psi_{23}$ is isotopic to the inclusion $M\hookrightarrow M'$.

Furthermore, we have
$$C_{h'}\circ C_{h}=\underline{\rm SHM}(\psi_{23}):\underline{\rm SHM}(-M,-\ga)\ra\underline{\rm SHM}(-M_2,-\ga_2).$$
\elem

\bpf
The two handles can be canceled topologically so the map $\psi_{23}$ is easy to find.

Suppose $(M_0,\ga_0)=(M,\ga)$. As in the definition \ref{defn_two_handle_gluing_map2} and definition \ref{defn_three_handle_gluing_map2}, for $i=0,1,2$, there are suitable closures $D_i=(Y_i,R_i,r_i,m_i,\eta_i)$ for $(M_i,\ga_i)$. The map $C_{h}$ is induced by a cobordism $-\widehat{W}$ from $-Y_0$ to $-Y_1$ so that $\widehat{W}$ is obtained from $Y_0\times[0,1]$ by attaching a $4$-dimensional $2$-handle $D^3\times[0,1]$  along a curve $\be\subset Y_0\times\{1\}$. The map $C_{h'}$ is induced by a cobordism $-\widehat{W}'$ from $-Y_1$ to $-Y_2$ so that $\widehat{W}'$ is obtained from $Y_1\times[0,1]$ by attaching a $4$-dimensional $3$-handle $D^3{}'\times[0,1]$ along a sphere in $Y_1\times\{1\}$. The $3$-dimensional handles $D^3$ and $D^3{}'$ is a pair of handles which can be canceled topologically, so the corresponding pair of $4$-dimensional handles will also be canceled topologically. Hence the composition of the cobordism $\widehat{W}\cup\widehat{W}'$ is actually diffeomorphic to the product cobordism $Y_0\times[0,2]$. We can think of identification $Y_0\times\{2\}$ with $Y_2$ to be induced by the diffeomorphism $\psi_{23}$ and hence by lemma \ref{lem_compare_canonical_maps} we have an equality
$$C_{h'}\circ C_{h}=\underline{\rm SHM}(\psi_{23})$$
\epf

We have a few more invariant results for contact handle attachments.

\blem\label{lem_independent_of_order}
Suppose $(M,\ga)$ is a balanced sutured manifold and $h=(\phi,S,D^3,\delta)$, $h'=(\phi',S',D^3{}',\delta')$ are two contact handles glued to $(M,\ga)$ with result $(M',\ga')$. Suppose further that the gluing maps have disjoint images:
$$\phi(S)\cap\phi'(S')=\emptyset,$$
then the two maps commute:
$$C_{h}\circ C_{h'}= C_{h'}\circ C_{h}:\underline{\rm SHM}(-M,-\ga)\ra\underline{\rm SHM}(-M',-\ga').$$
\elem

\bpf
Note that from definition \ref{defn_zero_and_one_handle_gluing_map}, \ref{defn_two_handle_gluing_map2} and \ref{defn_three_handle_gluing_map2}, we know that the gluing maps on canonical modules are essentially induced by cobordism which is either a product one or one obtained from a product cobordism by adding a $4$-dimensional $1$- or $2$-handle. The condition in the lemma means that the attachments of those $4$-dimensional handles can be moved apart and hence commute with each other. Hence the inducing gluing maps between canonical modules also commute.
\epf

\brem\label{rem_move_apart}
Suppose we first glue $h$ and then glue $h'$ so that the index of $h$ is no smaller than that of $h'$, then by an isotopy we call always move them apart. Hence such gluings always commute.
\erem

\blem\label{lem_commute_with_diffeomorphisms}
Suppose $(M,\ga)$ is a balanced sutured manifold with a local contact structure defined in a collar of $\partial{M}$. Suppose $(M',\ga')$ is another balanced sutured manifold. Suppose 
$$f:(M,\ga)\ra(M',\ga')$$
is a diffeomorphism and $f^{-1}$ will pull back the local contact structure. Suppose $h'$ is a contact handle attached to $(M',\ga')$ and via $f$ we can regard $h'$ as a contact handle $h$ attached to $(M,\ga)$ and there is a contactomorphism
$$\tilde{f}:M\cup h\ra M'\cup h'$$
which restricts to $f$ on $M$. Then we have an equality:
$$C_{h'}\circ\underline{\rm SHM}(f)=\underline{\rm SHM}(\tilde{f})\circ C_{h}.$$
\elem

\bpf
The instanton version is proved in \cite{baldwin2016instanton}. The monopole version is the same.
\epf

The last invariance result is about the inclusion of $(M,\ga)$ into a disjoint union $(M,\ga)\sqcup (N,\delta)$ when $(N,\delta)$ has a contact structure $\xi$ so that $\partial{N}$ is a convex surface and $\delta$ is the dividing set. Then from \cite{ozbagci2011contact} we know that $(N,\xi)$ possesses a contact handle decomposition $h_1,...,h_n$. We can regard those contact handles as attached to $(M,\ga)$ but all attaching maps are disjoint from $\partial{M}$. From this point of view, there is a map:
$$C_{h_n}\circ C_{h_{n-1}}\circ...\circ C_{h_1}:\underline{\rm SHM}(-M,-\ga)\ra\underline{\rm SHM}(-(M\sqcup N),-(\ga\cup\delta)).$$
We want to prove that this map is independent of the contact handle decomposition of $(N,\delta)$. The idea is that essentially this map is the identity on $M$ tensoring with the contact element of $N$. The proof will become easier if we require that there is no $3$-handle existing in the handle decomposition of $(N,\delta)$ as such decompositions can be related to partial open book decompositions of $(N,\delta)$. We will not introduce the basic definitions of partial open book decomposition or positive stabilizations, and interested readers are referred to \cite{baldwin2016contact,baldwin2016instanton,juhasz1803contact}.

\blem[Juh\'asz, Zemke, \cite{juhasz1803contact}, section 4.1]\label{lem_relation_between_open_book_and_handle_decomposition}
Suppose $(N,\delta)$ is a balanced sutured manifold and $\xi$ is a contact structure on $N$ so that $\partial{N}$ is convex and $\delta$ is the dividing set. The the following two objects are in one-to-one correspondence to each other:

(1). A partial open book decomposition of $(N,\delta,\xi)$.

(2). A handle decomposition of $(N,\delta,\xi)$ with no $3$-handles.
\elem

\blem[Honda, Kazez, Mati\'c, \cite{honda2009contact}, theorem 1.3]\label{lem_positive_stabilization}
Suppose $(M,\ga)$ is a balanced sutured manifold and $\xi$ is a positive contact structure on $(M,\ga)$ so that $\partial{M}$ is a convex surface and $\ga$ is the dividing set. Then $(M,\ga)$ admits a partial open book decomposition. Furthermore, for any two partial open book decompositions of $(M,\ga)$, one can perform positive stabilizations on each finitely many times so that the resulting two partial open book decompositions are isotopic.
\elem

\blem[Juh\'asz, Zemke, \cite{juhasz1803contact}, lemma 4.7]\label{lem_how_positive_stabilization_changes_handles}
Suppose $(N,\delta)$ is a balanced sutured manifold and $\xi$ is a contact structure on $N$ so that $\partial{N}$ is convex and $\delta$ is the dividing set. Suppose $(S,P,h)$ is a partial open book decomposition of $(N,\delta,\xi)$ and $(S',P',h')$ is a positive stabilization of $(S,P,h)$. Suppose $\mathcal{H}$ and $\mathcal{H}'$ are two contact handle decompositions of $(N,\delta,\xi)$ arising from $(S,P,h)$ and $(S',P',h')$ respectively, then $\mathcal{H}'$ can be obtained from $\mathcal{H}$ by adding a pair of canceling index $1$- and $2$-handles (See lemma \ref{lem_one_two_handle_cancelation}). 
\elem

\blem\label{lem_gluing_disjoint_union}
Suppose $(M,\ga)$ is a balanced sutured manifold and $(N,\delta)$ is a balanced sutured manifold with a compatible contact structure $\xi$. Suppose we have two different ways to decompose $(N,\delta,\xi)$ into contact handles: $h_1,...,h_n$ and $h_1',...,h_m'$, so that neither contains a $3$-handle. Then we can regard those handles as attached to $M$ and have an equality
$$C_{h_n}\circ...\circ C_{h_1}=C_{h_m'}\circ...\circ C_{h_1'}:\underline{\rm SHM}(-M,-\ga)\ra\underline{\rm SHM}(-(M\sqcup N),-(\ga\cup\delta)).$$
\elem

\bpf
The proof is a combination of lemmas \ref{lem_relation_between_open_book_and_handle_decomposition}, \ref{lem_positive_stabilization}, \ref{lem_how_positive_stabilization_changes_handles}, \ref{lem_one_two_handle_cancelation} and \ref{lem_commute_with_diffeomorphisms}.
\epf

\brem
This is essentially the way Baldwin and Sivek defined a contact invariant for sutured instantons in \cite{baldwin2016instanton}.
\erem

If we allow $3$-handles in the contact handle decompositions, the gluing map is still independent of decompositions and we will prove this result in the next proposition. Also we shall remark that this proof does not depend on the uniqueness part of the relative Giroux correspondence, but as a price to pay, local coefficients are necessary. Note that if one is already satisfied with using the uniqueness part of relative Giroux correspondence then the next technical proposition is not used anywhere else in the paper.

\bprop\label{prop_gluing_disjoint_union}
In lemma \ref{lem_gluing_disjoint_union}, if we allow $3$-handles in both of the decompositions, then the same conclusion still holds.
\eprop
\bpf
Since $\{h_i\}$ and $\{h_j'\}$ are both contact handle decompositions of $(N,\delta,\xi)$, they must both have at least one $0$-handle. Let $(N_0,\delta_0)$ be a $0$-handle or a $3$-ball with one simple closed curve being the suture on its boundary, we can view all other handles $h_i$ or $h_j'$ as being attached to $N_0$.

By lemma \ref{lem_independent_of_order} we can assume that the handles are ordered so that the index is non-decreasing. Suppose $(N_1,\ga_1)$ is gotten from $(N_0,\delta_0)$ by attaching all $0$- and $1$-handles in $\{h_1\}$ and $(N_2,\delta_2)$ is got from attaching all remaining $2$ and $3$-handles to $(N_1,\delta_1)$. There is a contactomorphism $g:N_2\ra N$. As in definition \ref{defn_marked_closure} and \ref{defn_closure}, we can use auxiliary data $T,f$ to form a pre-closure $\tilde{N}_1$ for $(N_1,\delta_1)$. We shall require:

(1). There is an arc configuration $\mathcal{A}$ (defined in \cite{baldwin2016contact}) on $T$ so that ${\widetilde{N}}_1$ carries a contact structure and $\partial{\widetilde{N}}_1$ is convex.

This can be achieved by an auxiliary surface $T$ of large enough genus. By definition \ref{defn_zero_and_one_handle_gluing_map}, we know that $\widetilde{N}_1$ is also a pre-closure for $(N_0,\delta_0)$. Define
$$\partial{\widetilde{N}_1}=R_+\cup R_-.$$

From definition \ref{defn_two_handle_gluing_map2} and \ref{defn_three_handle_gluing_map2}, there are curves $\be_1,...,\be_s\subset {\rm int}(N_1)\subset \widetilde{N}_1$, all isotopic to curves on $\partial{N_1}$, and spheres $S_1,...,S_u\subset {\rm int}(N_1)\subset Y_1$, so that if we do $0$-surgeries along all $\be_i$ and do cut and paste surgeries along all $S_j$, then we will get a pre-closure $\widetilde{N}_2$ for $(N_2,\delta_2)$ so that $\partial{\widetilde{N}_2}=\partial{\widetilde{N}_1}$. As discussed in proposition \ref{prop_equivalence_of_the_definition}, the Dehn surgeries can be made to be contact $+1$ surgeries and the cut and paste surgeries can also be done to preserve contact structures. So there will be a contact structure on $\widetilde{N}_2$.

We can similarly form $(N_1',\delta_1')$, $(N_2',\delta_2')$ and a contactomorphism $g':N_2'\ra N$. We can require

(2). The pre-closure $\widetilde{N}_1$ is also a pre-closure for $(N_1',\delta_1')$.

This requirement can be achieved by choosing a $T$ with large enough genus. 

As above there are curves $\be_i'\subset \widetilde{N}_1$ and spheres $S_j'\subset\widetilde{N}_1$ so that doing suitable surgeries along these objects will result in a pre-closure $\widetilde{N}_2'$ carrying a suitable contact structure and $\partial{\widetilde{N}_2}'=\partial{\widetilde{N}_1}$. Note the boundary of all pre-closures are identified and are all $R_+\cup R_-$. For later use, we will need a diffeomorphism
$$C:\tilde{N}_2\ra \widetilde{N}_2'.$$
We shall require that 
$$C|_{(N_2)}=(g')^{-1}\circ g:N_2\ra N_2'.$$

If we pick a non-separating simple closed curve $c\subset T$, then $c$ will correspond to two curves $c_+\subset R_+$ and $c_-\subset R_-$. By choosing the auxiliary surface $T$ with large enough genus, we could require that 

(3). $C$ preserves $c\times[-1,1]\subset \widetilde{N}_2$.

(4). There exists a gluing diffeomorphism $h:R_+\ra R_-$ preserving contact structures and identifying $c_+$ with $c_-$ in the way that $c_+$ and $c_-$ are both identified with $c$.

(5). There exists a smooth curve $\eta\subset R$ intersecting $c$ transversely once.

We can use the same auxiliary data $R_+,h,\eta$ to get a marked closure $\mathcal{D}_0=(Y_0,R_+,r,m_0,\eta)$ for $(N_0,\delta_0)$ and marked closures $\mathcal{D}_2=(Y_2,R_+,r,g^{-1},\eta)$, $\mathcal{D}_2'=(Y_2',R_+,r,(g')^{-1},\eta)$ for $N$.

Note that for all three marked closures, the curve $c$ becomes a torus as follows: there are $c\times [-1,1]\subset T\times [-1,1]$ and $c\times[-1,1]\in R_+\times[-1,1]$ and the way we get the closures will identify their boundaries to get $c\times S^1$. Suppose $\Sigma_0\in Y_0$, $\Sigma_2\subset Y_2$ and $\Sigma_2'\subset Y_2'$ are the corresponding tori.

 Now let $\widetilde{T}$ be an auxiliary surface for $(M,\ga)$ and $\tilde{f}:\partial{\widetilde{T}}\ra\ga$ be the gluing map. Let
$$\widetilde{M}=M\mathop{\cup}_{\tilde{f}}\widetilde{T}\times[-1,1],~\partial{\widetilde{M}}=\widetilde{R}'_+\cup \widetilde{R}_-.$$
Suppose there is a non-separating simple closed curve $\tilde{c}\subset \widetilde{T}$ and a diffeomorphism
$$\tilde{h}: \widetilde{R}_+\ra \widetilde{R}_-$$
so that $\tilde{h}(\tilde{c}\times \{1\})=\tilde{c}\times \{-1\}$. Use $\tilde{h}$ we will get a marked closure
$$\mathcal{D}=(Y,\tilde{R}_+,r,m,\tilde{\eta})$$
for $(M,\ga)$ and there is a torus $\Sigma\subset Y$ corresponding to $\tilde{c}$. 

We can now form a marked closure for $(M,\ga)\sqcup (N_0,\delta_0)$ as follows. Cut $Y$ open along $\Sigma$, and let $Y''=Y\backslash ({\rm int}(N(\Sigma))).$ We have $\partial Y''=\Sigma_+\cup\Sigma_-$. Cut $Y_0$ along $\Sigma_0$ and let $Y_0''=Y_0\backslash ({\rm int}(N(\Sigma_0)))$ with $\partial Y_0''=\Sigma_{0,+}\cup\Sigma_{0,-}$. Let
$$\tau:\Sigma\ra\Sigma_0$$
be a diffeomorphism so that $\tau(\Sigma\cap\tilde{\eta})=\Sigma_0\cap\eta.$ We can use $\tau$ to glue $\Sigma_+$ to $\Sigma_{0,-}$ and $\Sigma_-$ to $\Sigma_{0,+}$. Let $\hat{Y}$ be the resulting manifold. There are corresponding $\widehat{R},\hat{r},\hat{m},\hat{\eta}$ so that
$$\widehat{\mathcal{D}}=(\widehat{Y},\widehat{R},\hat{r},\hat{m},\hat{\eta})$$
is a marked closure of $(M,\ga)\sqcup (N_0,\ga_0)$. If we use $Y_2$ or $Y_2'$ and the same $\tau$, we can construct two similar marked closures
$$\widehat{\mathcal{D}}_2=(\widehat{Y}_2,\widehat{R},\hat{r},\hat{m}_2,\hat{\eta}),\widehat{\mathcal{D}}_2'=(\widehat{Y}_2',\widehat{R},\hat{r},\hat{m}_2',\hat{\eta})$$
for $(M,\ga)\sqcup (N,\delta)$. Now the diffeomorphism $C$ extends by identity to a diffeomorphism which we also called $C$:
$$C:\widehat{Y}_2\backslash{\rm int}({\rm im}(\hat{r}))\ra\widehat{Y}_2'\backslash{\rm int}({\rm im}(\hat{r}')).$$ 

There are Legendrian curves and spheres $\be_i,S_j\subset (Y\backslash N(\Sigma))\subset \widehat{Y}$ so that if we do contact $+1$-surgeries along these curves $\be_i$ and do cut and paste surgeries along $S_j$, then the resulting manifold will be exactly $\widehat{Y}_2$. Hence there is a cobordism $W$ from $-\widehat{Y}$ to $-\widehat{Y}_2$ so that $W$ is obtained obtained from $\widehat{Y}\times[0,1]$ by gluing $0$-framed $4$-dimensional $2$-handles along all $\be_i\subset\widehat{Y}\times\{1\}$, and gluing $4$-dimensional $3$-handles along all $S_j\subset \widehat{Y}\times \{1\}$. Then the map
$$HM(W):HM(-\widehat{Y}|-\hat{r}(\{0\}\times \widehat{R}),\Gamma_{-\hat{\eta}})\ra HM(-\widehat{Y}_2|-\hat{r}(\{0\}\times \widehat{R}),\Gamma_{-\hat{\eta}})$$
will induce the map
$$C_{h_n}\circ...\circ C_{h_2}:\underline{\rm SHM}(-(M\sqcup N_0),-(\ga\cup\delta_0))\ra\underline{\rm SHM}(-(M\sqcup N),-(\ga\cup\delta)).$$

Later we will use another interpretation of $-W$. Gluing $4$-dimensional $2$- and $3$- handles to $\widehat{Y}\times[0,1]$ at $\widehat{Y}\times\{1\}$ is equivalent to glue $4$-dimensional $2$- and $1$-handles to $\widehat{Y}_2\times [0,1]$ at $\widehat{Y}_2\times\{0\}$. Suppose those handles are attached along curves $\theta_i$, which correspond to $\be_i$, and along pairs of points $(p_j,q_j)$, which correspond to $S_j$.

There are curves $\be_j'\subset \widehat{Y}$ and spheres $S_j'\subset\widehat{Y}'_2$ as well. We can construct similarly a cobordism $W'$ from $-\widehat{Y}$ to $-\widehat{Y}_2'$ which induces the map
$$C_{h_m'}\circ...\circ C_{h_2'}:\underline{\rm SHM}(-(M\sqcup N_0),-(\ga\cup\delta_0))\ra\underline{\rm SHM}(-(M\sqcup N),-(\ga\cup\delta)).$$

Just as for $W$, there are curves $\theta_i'\subset \widehat{Y}_2'$ corresponding to $\be_i'$ and pairs of points $(p_j',q_j')$ corresponding to $S_j'\subset \widehat{Y}_2'$.

To show that
$$C_{h_n}\circ...\circ C_{h_2}=C_{h_m'}\circ...\circ C_{h_2'},$$
we only need to show that
\begin{equation}\label{eq_1}
\Phi^g_{-\widehat{\mathcal{D}}_2,-\widehat{\mathcal{D}}_2'}\circ HM(W)\doteq HM(W'),
\end{equation}
where $\Phi_{-\widehat{\mathcal{D}}_2,-\widehat{\mathcal{D}}'_2}$ is the canonical map constructed in definition \ref{defn_caninical_map_for_the_same_genus} for the two marked closures $-\widehat{\mathcal{D}_2}$ and $-\widehat{\mathcal{D}}'_2$ of the same genus. The diffeomorphism $C$ is used to construct such a canonical map. As in definition \ref{defn_caninical_map_for_the_same_genus}, we can construct
$$\varphi^C:\widehat{R}\ra\widehat{R},~\varphi^C:\widehat{R}\ra\widehat{R}.$$
The way we choose $C$ makes sure that the maps $\varphi^C$ and $\psi^C$ will fix the part of $\widehat{R}$ coming from $\widetilde{R}_+$ which was used to build the marked closure $\mathcal{D}$ of $(M,\ga)$. Hence we can decompose $\varphi^C\circ\psi^C$ and $(\psi^C)^{-1}$ as 
$$\varphi^C\circ\psi^C\sim D^{e_1}_{a_1}\circ...\circ D^{e_p}_{a_p},~\varphi^C\circ\psi^C\sim D^{e_{u+1}}_{a_{u+1}}\circ...\circ D^{e_q}_{a_q},$$
so that all $a_k$ are disjoint from the part of $\widehat{R}$ coming from  $\widetilde{R}_+$. 

In general there will be both positive and negative Dehn twists but for simplicity, we only deal with the case when all $e_i=-1$. The general case will follow from a similar argument. As in definition \ref{defn_caninical_map_for_the_same_genus}, let $W_{c}$ be the cobordism from $-\widehat{Y}_2$ to $-\widehat{Y}_2'$ obtained from $(-\widehat{Y}_2)\times [0,1]$ by gluing some $4$-dimensional $2$-handles along all the curves $a_1,...,a_q\subset (-\widehat{Y}_2)\times \{1\}$, then the canonical map $\Phi_{-\widehat{\mathcal{D}}_2,-\widehat{\mathcal{D}}_2'}$ is induced by the cobordism $W_c$.

Equation (\ref{eq_1}) is now equivalent to
\begin{equation}\label{eq_2}
	HM(W\cup W_c)\doteq HM(W').
\end{equation}

Note that the curves $\theta_i$ and pairs of points $\{p_j,q_j\}$ used to define $W$ are all contained in $Y_0\backslash\mathring{N}(\Sigma_0)\subset \hat{Y}$, so intuitively there is nothing happened in $Y$ part of $\hat{Y}$ and we shall be able to split off a product copy of $Y$. This idea is carried out explicitly as follows.

\begin{figure}[h]
\centering
\begin{overpic}[width=5.0in]{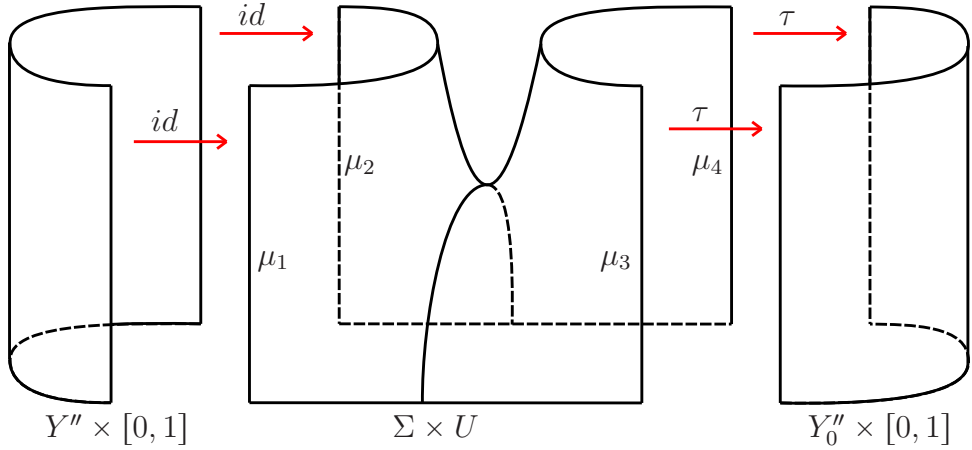}
	\put(4,-3){$Y''\times[0,1]$}
	\put(40,-3){$\Sigma\times U$}
	\put(83,-3){$Y_0''\times[0,1]$}
	\put(15,29){$id$}
	\put(24,40){$id$}
	\put(71,30){$\tau$}
	\put(80,40){$\tau$}
	\put(35,25){$\mu_2$}
	\put(26,15){$\mu_1$}
	\put(71,25){$\mu_4$}
	\put(61.5,15){$\mu_3$}
\end{overpic}
\vspace{0.05in}
\caption{The three parts of the cobordism $-W_e$. The middle part is $\Sigma\times U$, while the $\Sigma$ directions shrink to a point in the figure.}\label{fig_excision_cobordism}
\end{figure}

Let $U$ be the surface depicted as in the figure \ref{fig_excision_cobordism}. It has four vertical  parts of the boundary which we call $\mu_1,...,\mu_4$. Suppose each is parametrized by $[0,1]$. Recall we have $Y''$ and $Y_0''$ by cutting open $Y$ along $\Sigma$ and $Y_0$ along $\Sigma_0$ repsectively. We have the gluing diffeomorphism $\tau$ to get $\widehat{Y}$. Now let $-W_e$ be the cobordism obtained by gluing three parts $Y''\times[0,1]$, $\Sigma\times U$ and $Y_0''\times [0,1]$ where we use $id\times id$ to glue $\partial{Y''}\times [0,1]$ to $\Sigma\times (\mu_1\cup\mu_2)$ and use $\tau\times id$ to glue $\Sigma\times (\mu_3\cup\mu_4)$ to $\partial{Y_0''}\times [0,1]$. The result $W_e$ can be thought of as a cobordism from $-(Y\sqcup Y_0)$ to $-\widehat{Y}$. Similarly we can construct a cobordism $W_e'$ from $(Y\sqcup Y_2')$ to $\widehat{Y}'_2$. The same cobordism is juts one from $-\widehat{Y}'_2$ to $-(Y\sqcup Y_2')$. From theorem 3.2 in \cite{kronheimer2010knots}, we know that $W_e$ and $W_e'$ induces isomorphisms, so the equality (\ref{eq_2}) is equivalent to
\begin{equation}\label{eq_3}
	HM(W_e\cup W\cup W_c\cup W_e')\doteq HM(W_e\cup W'\cup W_e').
\end{equation}

\begin{figure}[h]
\centering
\begin{overpic}[width=5.0in]{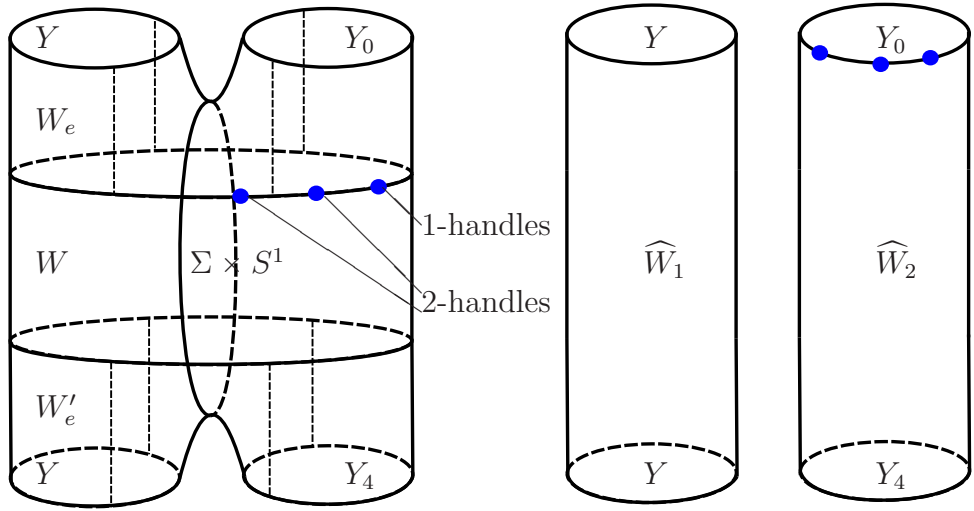}
	\put(3,3){$Y$}
	\put(3,10){$W_e'$}
	\put(3,25){$W$}
	\put(3,40){$W_e$}
	\put(3,48.5){$Y$}
	\put(35,3){$Y_4$}
	\put(35,48.5){$Y_0$}
	\put(19,25){$\Sigma\times S^1$}
	\put(39,34){\line(1,-1){4}}
	\put(43,29){$1$-handles}
	\put(33,33){\line(1,-1){10}}
	\put(43,21){$2$-handles}
	\put(25,33){\line(3,-2){18}}
	\put(66,3){$Y$}
	\put(66,48.5){$Y$}
	\put(90,3){$Y_4$}
	\put(90,48.5){$Y_0$}
	\put(66,25){$\widehat{W}_1$}
	\put(90,25){$\widehat{W}_2$}
\end{overpic}
\vspace{0.05in}
\caption{Cut along $\Sigma\times S^1$ and glue back two copies of $\Sigma\times D^2$.}\label{fig_split_apart}
\end{figure}

On $W_e\cup W\cup W_c\cup W_e'$ we can cut along a $3$-manifold $\Sigma\times S^1$ as shown in the figure \ref{fig_split_apart}, and glue back two copies of $\Sigma\times D^2$ along boundaries. The result is a cobordism $\widehat{W}=\widehat{W}_1\cup\hat{W}_2$, where $\hat{W}_1\cong(-Y)\times[0,1]$, and $\hat{W}_2$ is a cobordism from $-Y_0$ to $-Y_2'$ obtained from $Y_2'\times [0,1]$ by attaching $4$-dimensional $1$-handles at pairs of points $\{p_j,q_j\}\subset Y_2\times\{1\}$ and then attaching $4$-dimensional $2$-handles to $Y_2\times[0,1]$ along curves $\theta_i\subset Y\times\{1\}$ and $a_k\subset Y_2\times\{1\}$. We can apply similar argument to $W_e\cup W'\cup W_e'$, and get $\widehat{W}'=\widehat{W}_1'\cup \widehat{W}_2'$, where $\widehat{W}_1'=\widehat{W}_1\cong(-Y)\times [0,1]$ and $\hat{W}_2'$ is a cobordism from $-Y_0$ to $-Y_2'$ obtained from $Y_2'\times [0,1]$ by attaching $4$-dimensional $1$-handles at $\{p_j',q_j'\}\subset Y_2'\times\{1\}$ and then attaching $4$-dimensional $2$-handles to $Y_2'\times[0,1]$ along $\theta_i'\subset Y_2\times\{1\}$. There is a commutative diagram from the naturality of K\"unneth formula: 
\begin{equation*}
\xymatrix{
HM(-(Y\sqcup Y_0))\ar[dd]^{HM(\widehat{W})'}\ar[dd]_{HM(\widehat{W})}\ar[rr]^{i}&&HM(-Y)\otimes HM(-Y_0)\ar[dd]_{HM(\widehat{W}_1)\otimes HM(\widehat{W}_2)}\ar[dd]^{HM(\widehat{W}_1')\otimes HM(\widehat{W}_2')}\\
&&\\
HM(-(Y\sqcup Y_2'))\ar[rr]^{i'}&&HM(-Y)\otimes HM(-Y_2')\\
}	
\end{equation*}

The map $HM(\widehat{W}_1)=HM(\widehat{W}_1')$ since they are both product cobordisms. We claim that $HM(\widehat{W}_2)\doteq HM(\widehat{W}_2')$. Since both cobordisms are exact symplectic, they both map contact elements to contact elements (see corollary 2.23 in \cite{baldwin2016contact}). Yet $Y_0$ is a surface fibration over $S^1$ with fibre $R_0$, hence $HM(Y_0)\cong \mathcal{R}$ (see lemma 4.7 and 4.9 in \cite{kronheimer2010knots}), and the contact element in $HM(Y_0)$ is a generator of the module (See \cite{baldwin2016contact}). From K\"unneth formula, the maps $i$ and $i'$ are injective, so $HM(\widehat{W})\doteq HM(\widehat{W}')$. Finally from corollary 2.10 (or see the proof of theorem 3.2) of \cite{kronheimer2010knots}, we know that 
$$HM(W_e\cup W\cup W_c\cup W_e')\doteq HM(\widehat{W})\doteq HM(\widehat{W}')\doteq HM(W_e\cup W'\cup W_e').$$
Hence we are done.
\epf

\section{The general gluing maps}
Now we will try to construct the general gluing map.

\bdefn\label{defn_sutured_submanifold}
Suppose $(M',\ga')$ is a balanced sutured manifold. By {\it sutured submanifold} we mean a balanced sutured manifold $(M,\ga)$ so that $M\subset {\rm int}(M')$.
\edefn

In \cite{juhasz1803contact} Juh\'asz and Zemke used contact cell decompositions to re-construct the gluing map originally introduced by Honda, Kazez and Mati\'c in \cite{honda2008contact}. Here we will introduce the basic definition of contact cell decompositions and use it to construct general gluing maps. The following definition is from \cite{juhasz1803contact}.

\bdefn\label{defn_contact_cell_decomposition}
Suppose $(M,\ga)$ is a sutured submanifold of $(M',\ga')$ and $\xi$ is a contact structure on $(Z=M\backslash{\rm int}(M),\ga\cup\ga')$, so that $\partial{Z}$ is convex and $\ga\cup\ga'$ is the dividing set. A {\it contact cell decomposition} of $(Z,\xi)$ consists of the following data:

(1). A non-vanishing contact vector field $v$ that is defined on a neighborhood of $\partial{Z}\subset Z$ and with respect to which $\partial{Z}$ is a convex surface with dividing set $\ga\cup\ga'$. 

The flow of $v$ induces a diffeomorphism between $\partial{Z}\times I$ and a collar neighborhood of $\partial{Z}$ and under this diffeomorphism $\nu$ corresponds to the vector field $\frac{\partial}{\partial t}$, $\partial{M}$ is identified with $\partial{M}\times\{0\}$ and $\partial{M'}$ is identified with $\partial{M}'\times\{1\}.$ We shall call
$$\nu=v|_{\partial{M}\times I},~\nu'=v|_{\partial{M'}\times I}.$$

(2). Barrier surfaces 
$$S\subset \partial{M}\times (0,1),~S'\subset \partial{M'}\times (0,1)$$
that are isotopic to $\partial{M},\partial{M'}$ respectively and are transverse to $v$. Write $N$ for the collar neighborhood of $\partial{M}$ bounded by $S$ and $N'$ for $\partial{M}',S'$ similarly. We shall call
$$Z'=Z\backslash{\rm int}(N\cup N').$$ Note $\partial{Z}'=S\cup S'$ is a convex surface. 

(3). A Legendrian graph $\Gamma\subset Z'$ which intersects $\partial{Z'}$ transversely in a finite collection of points along the dividing set on $\partial{Z}'$ with respect to $v$. Furthermore, $\Gamma$ is tangent to $v$ in a neighborhood of  $\partial{Z}'\subset Z'$.

(4). A choice of regular neighborhood $N(\Gamma)\subset Z'$ of $\Gamma$ such that $\xi$ is tight on $N(\Gamma)$ and $\partial{N(\Gamma)}\backslash \partial{Z}'$ is a convex surface. We also require that $N(\Gamma)\cap\partial{Z'}$ is a collection of disks $D$ with Legendrian boundary such that each boundary $\partial{D}$ has $tb(D)=-1$. We shall also assume that $N(\Gamma)$ meets $\partial{Z}'$ tangentially along the Legendrian unknots forming $\partial{N(\Gamma)}$.

(5). A collection of $2$-cells $D_1,...,D_n$ inside $Z'\backslash{\rm int}(N(\Gamma))$ with Legendrian boundary on $\partial{Z}'\cup\partial{N(\Gamma)}$ and each $\partial{D}_i$ has $tb(\partial{D}_i)=-1$. Furthermore, the following two conditions shall hold:

(a). Each component of
$Z'\backslash(N(\Gamma)\cup D_1\cup...\cup D_n)$
is a $3-$ball and $\xi$ is tight on each of them.

(b). The disks $N(\Gamma)\cap \partial{Z}$ and the Legendrian arcs $\partial{D_i}\cap \partial{Z}'$ induces a sutured cell decomposition (\cite{juhasz1803contact}, definition 3.1), with the dividing set induced by $v$.
\edefn

Now we are ready to define the gluing map:
\bdefn\label{defn_contact_gluing_map_for_monopoles}
Suppose $(M,\ga)$ is a sutured submanifold of a balanced sutured manifold $(M',\ga')$. Suppose $Z=M'\backslash{\rm int}(M)$ is equipped with a contact structure $\xi$ so that $\partial{Z}$ is a convex surface with dividing set $\ga\cup\ga'$. Suppose $\mathcal{C}$ is a contact cell decomposition of $(Z,\xi)$, we will use the same notations as in definition \ref{defn_contact_cell_decomposition}. The contact vector field $\nu$ will induce a diffeomorphism $\phi_{\nu}: (M,\ga)\ra(M\cup N,\delta)$ where $\delta$ is the dividing set on $S\subset \partial{N}$ with respect to $\nu$. Suppose $h_1',...,h_n'$ is a handle decomposition of $(N',\delta'\cup\ga')$ with no $3$-cells. The existence of such decomposition is guaranteed by lemma \ref{lem_positive_stabilization} and lemma \ref{lem_relation_between_open_book_and_handle_decomposition}.

A contact cell decomposition will lead to a handle decomposition: vertices of $\Gamma$ are $0$-handles, edges of $\Gamma$ are $1$-handles, $2$-cells $D_i$ are $2$-handles and the remaining is a collection of $3$-handles. Suppose we get a sequence of contact handles ${h}'_1,...,{h}'_m$ from it, then we define the contact gluing map
$$\Phi_{\xi}:\underline{SHM}(-M,-\ga)\ra \underline{SHM}(-M',-\ga')$$
to be 
$$\Phi_{\xi}= C_{{h}'_m}\circ...\circ C_{{h}'_1}\circ C_{h_n}\circ...\circ C_{h_1}\circ\underline{SHM}(\phi_{\nu}).$$
\edefn

\bprop\label{prop_well_definedness_of_gluing_map}
The contact gluing map $\Phi_{\xi}$ as in definition \ref{defn_contact_gluing_map_for_monopoles} is well defined.
\eprop

\bpf
The relation between two contact cell decompositions is stated in \cite{juhasz1803contact}, proposition 3.6. Any two contact cell decompositions are actually related by a sequence of isotopies fixing boundary and three types of cancelations. The well-definedness of our gluing map is just a combination of that proposition with lemmas \ref{lem_zero_one_cancelation}, \ref{lem_one_two_handle_cancelation}, \ref{lem_two_three_handle_cancelation}, \ref{lem_commute_with_diffeomorphisms} and \ref{lem_gluing_disjoint_union}.  
\epf

\bprop\label{prop_functoriality_of_gluing_map}
Suppose $(M,\ga)$ is a sutured submanifold of $(M',\ga')$ and $(M',\ga')$ is a sutured submanifold of $(M'',\ga'')$. Suppose there are contact structures $\xi$ on $Z=M'\backslash{\rm int}(M)$ and $\xi'$ on $Z'=M''\backslash{\rm int}(M')$, and their union $\xi''=\xi\cup\xi'$ is a contact structure on $Z''=M''\backslash{\rm int}(M)$, so that the boundaries of corresponding manifolds are all convex surfaces and the sutures are dividing sets. Then we have an equality:
$$\Phi_{\xi'}\circ\Phi_{\xi}=\Phi_{\xi''}:\underline{SHM}(-M,-\ga)\ra \underline{SHM}(-M'',-\ga'').$$
\eprop

\bpf
We follow the idea from Juh\'asz and Zemke \cite{juhasz1803contact}. Suppose $\mathcal{C}$ is a contact cell decomposition of $Z$ and $\nu$ is defined as in definition \ref{defn_contact_cell_decomposition}. Suppose $h_1,...,h_n$ and $h_1',...,h_m'$ are defined using $\mathcal{C}$ as in definition \ref{defn_contact_gluing_map_for_monopoles}. Suppose $\mathcal{C}'$ is a contact cell decomposition of $Z'$ with $\tilde{\nu}$, $\tilde{h}_1,...,\tilde{h}_s$ and $\tilde{h}'_1,...,\tilde{h}'_t$ defined similarly. Then we have
$$\Phi_{\xi'}\circ\Phi_{\xi}=C_{\tilde{h}'_t}\circ...\circ C_{\tilde{h}'_1}\circ C_{\tilde{h}_s}\circ...\circ C_{\tilde{h}_1}\circ\underline{SHM}(\phi_{\tilde{\nu}})\circ C_{{h}'_m}\circ...\circ C_{{h}'_1}\circ C_{h_n}\circ...\circ C_{h_1}\circ\underline{SHM}(\phi_{\nu}).$$

Suppose $\bar{h}_i=\phi_{\tilde{\nu}}(h_i)$ and $\bar{h}'_i=\phi_{\tilde{\nu}}(h'_i)$, then by lemma \ref{lem_commute_with_diffeomorphisms}, we know that
\begin{equation*}
	\underline{SHM}(\phi_{\tilde{\nu}})\circ C_{{h}'_m}\circ...\circ C_{{h}'_1}\circ C_{h_n}\circ...\circ C_{h_1}=C_{\bar{h}'_m}\circ...\circ C_{\bar{h}'_1}\circ C_{\bar{h}_n}\circ...\circ C_{\bar{h}_1}.
\end{equation*}

If we go back to the definition of gluing maps, we can see that the set of handles $\tilde{h}_1,...,\tilde{h}_s$ and the set of handles $\bar{h}_1,...,\bar{h}_n,\bar{h}_1',...,\bar{h}_m'$ are attached to disjoint parts, so we can switch their order by lemma \ref{lem_independent_of_order}. The handles $\tilde{h}_1,...,\tilde{h}_s$ corresponding to the neighborhood of $\partial{M}''$ bounded by $\partial{M}''$ and the barrier surface $\tilde{S}'\subset Z''$. The handles $\tilde{h}'_1,...,\tilde{h}'_t$ and $\bar{h}_1',...,\bar{h}_m'$ corresponding to Legendrian graphs and $2$-cells and tight $3$-balls in $Z$ and $Z'$, so they are still basic elements to form a contact cell decomposition of $Z\cup Z'$. The remaining handles $h_1,...,h_n$ correspond to the neighborhood of $\partial{M}'$ in $Z$ bounded by $\partial{M}'$ and the barrier surface $S'$. They consist of only $0$-, $1$- and $2$- handles by lemma \ref{lem_positive_stabilization} so we can consider them as Legendrian graphs and $2$-cells. Hence the whole series
$$C_{\tilde{h}'_t}\circ...\circ C_{\tilde{h}'_1}\circ C_{\tilde{h}_s}\circ...\circ C_{\tilde{h}_1}\circ C_{\bar{h}'_m}\circ...\circ C_{\bar{h}'_1}\circ C_{\bar{h}_n}\circ...\circ C_{\bar{h}_1}\circ\underline{SHM}(\phi_{\nu})$$
can be thought of as from some contact cell decomposition of $Z\cup Z'$ and hence the proposition follows.
\epf

Suppose $(M,\ga)$ is a sutured submanifold of $(M',\ga')$ and if $Z=M'\backslash{\rm int}(M)$ is just a product $\partial{M}\times[0,1]$ equipped with an $I$-invariant contact structure so that for any $t\in[0,1]$, $\partial{M}\times\{t\}$ is convex with $\ga\times\{t\}$ being the dividing set, Then we shall expect the contact gluing map to be the 'identity'. This is made precise by the following proposition.

\bprop\label{prop_trivial_gluing_map}
Suppose $(M,\ga)$ is a sutured submanifold of $(M',\ga')$ and $\xi$ is a compatible contact structure on $Z=M'\backslash{\rm int}(M)$. Suppose there is a Morse function $f$ and a contact vector field $nu$ on $Z$ so that

(1). There is no critical point of $f$ and 
$$f(\partial{M})=0,~f(\partial{M}')=1.$$

(2). The contact vector field $\nu$ is gradient like: $\nu(f)>0$ everywhere in $Z$.

Then we have the equality
$$\Phi_{\xi}=\underline{SHM}(\phi_{\nu}):\underline{SHM}(-M,-\ga)\ra\underline{SHM}(-M',-\ga').$$
where $\phi_{\nu}$ is just the diffeomorphism induced by $\nu$.
\eprop

\bpf
With lemma \ref{lem_zero_one_cancelation}, \ref{lem_one_two_handle_cancelation} and \ref{lem_two_three_handle_cancelation}, the proof is exactly the same as the proof of proposition 5.1 in \cite{juhasz1803contact}.
\epf

At the end of the section, we want to relate the general gluing map with the contact handle gluing map introduced before. Suppose $(M,\ga)$ is balanced sutured manifold and $h$ is a contact handle attached to $(M,\ga)$ and $(M',\ga')$ is the result balanced sutured manifold. First we shall note that $(M,\ga)$ is not a sutured submanifold as in definition \ref{defn_sutured_submanifold}, we require $M\subset{\rm int}(M')$. The way to resolve this is to glue a product region $\partial{M}\times [0,1]$ to $M$ along $\partial{M}\times\{0\}$ and glue $h$ to $\partial{M}\times \{1\}$. This is made precise by the following definition from \cite{juhasz1803contact}:
\bdefn\label{defn_morse_type_contact_handle}
Suppose $(M,\ga)$ is a sutured submanifold of $(M',\ga')$ and $\xi$ is a compatible contact structure on $Z=M'\backslash{\rm int}(M)$. Suppose there is a contact vector field $\nu$ on $Z$ and a decomposition $Z=Z_0\cup h$ such that

(1). The contact vector field $\nu$ points into $Z$ on $\partial{M}\subset \partial{Z}$ and points out of $Z$ on $\partial{M}'\subset \partial{Z}$.

(2). We have $Z_0\cong \partial{M}\times [0,1]$ and $\partial{M}$ is identified with $\partial{M}\times \{0\}\subset \partial{Z}$. We shall also require that $\nu$ is non-vanishing on $Z_0$, pointing into $Z_0$ on $\partial{M}\times \{0\}$, pointing out of $Z_0$ on $\partial{M}\times \{1\}$ and each flow line of $\nu$ on $Z_0$ is an arc from $\partial{M}\times \{0\}$ to $\partial{M}\times \{1\}$.

(3). We shall require that $h$ is a topologically $3$-ball with piece-wise smooth boundary and is tight under $\xi$.

(4). We can view $h$ as a contact $k$-handle, for $k=0,1,2,3$, attached to $M\cup Z_0$, with corner smoothed.

Then $(Z,\xi)$ is called a {\it Morse-type contact handle of index $k$}.
\edefn

\bprop\label{prop_equivalence_between_handle_and_general_gluing}
Suppose $(M,\ga)$ is a sutured submanifold of $(M',\ga')$ $(Z=M'\backslash{\rm int}(M),\xi)$ is a Morse-type contact handle of  index $k$ for $k=0,1,2,3$. Suppose the contact vector field $\nu$ and the decomposition $Z=Z_0\cup h$ are as in the definition \ref{defn_morse_type_contact_handle} and $\ga_0\subset\partial{M}\times \{1\}\subset\partial{Z_0}$ is the dividing set with respect to $\nu$. Suppose $\phi_{\nu}:(M,\ga)\ra (M\cup Z_0,\ga_0)$ is the diffeomorphism induced by $\nu$. Then we have an equality
$$\Phi_{\xi}=C_{h}\circ\underline{SHM}(\phi_{\nu}):\underline{SHM}(-M,-\ga)\ra\underline{SHM}(-M',-\ga').$$
\eprop
\bpf
The proof is exactly the same as the proof of proposition 5.6 in \cite{juhasz1803contact}. The handle cancelations needed have been proved in lemma \ref{lem_zero_one_cancelation}, \ref{lem_one_two_handle_cancelation} and \ref{lem_two_three_handle_cancelation}. 
\epf

\brem
With proposition \ref{prop_equivalence_between_handle_and_general_gluing}, \ref{prop_functoriality_of_gluing_map}, \ref{prop_well_definedness_of_gluing_map}, we can actually prove the conjecture 1.7 in \cite{baldwin2016contact}: the gluing maps constructed by composing contact handle gluing maps is independent of the contact handle decomposition.
\erem

\bcor
The contact element in sutured monopole Floer homology is preserved by the gluing map $\Phi_{\xi}$.
\ecor

\section{The cobordism maps}
\subsection{Constructions and functoriality}
Now we are ready to construct the cobordism map between sutured monopole Floer homologies. The following definitions are from \cite{juhasz2016cobordisms}.

\bdefn\label{defn_equivalent_contact_structures}
Suppose $(M,\ga)$ is a balanced sutured manifold and $\xi_0$, $\xi_1$ are two compatible contact structures. We say that $\xi_0$ and $\xi_1$ are {\it equivalent} if there is a $1$-parameter family $\xi_t$ so that for any $t\in[0,1]$, $\xi_t$ is a contact structure on $M$ with convex boundary $\partial{M}$.
\edefn

\bdefn\label{defn_suture_cobordism}
Suppose $(M_0,\ga_0)$ and $(M_1,\ga_1)$ are two balanced sutured manifolds. A {\it sutured cobordism} from $(M_0,\ga_0)$ to $(M_1,\ga_1)$ is a triple $\mathcal{W}=(W,Z,\xi)$ so that

(1). $W$ is a compact $4$-dimensional smooth oriented manifold with boundary

(2). $Z$ is a compact oriented $3$-manifold so that $\partial{W}\backslash {\rm int}(Z)=-M_1\cup M_2$.

(3). We have that $\xi$ is an oriented and co-oriented contact structure on $Z$ so that $\partial{Z}=\partial{M}_1\cup \partial{M_2}$ (not specifying the orientation) is a convex surface with dividing set $\ga_0\cup\ga_1$.
\edefn

\bdefn
Suppose $(M_0,\ga_0)$ and $(M_1,\ga_1)$ are two balanced sutured manifolds and $\mathcal{W}=(W,Z,\xi)$ is a suture cobordism between them. We can regard $(M_0,\ga_0)$ as a sutured submanifold of $(M_0\cup(-Z),\ga_1)$, and from definition \ref{defn_contact_gluing_map_for_monopoles} we have a gluing map
$$\Phi_{-\xi}:\underline{\rm SHM}(M_0,\ga_0)\ra\underline{\rm SHM}(M_0\cup(-Z),\ga_1).$$
The cobordism $W$ can be thought as one with sutured surface $(\partial{M_2},\ga_2)$, from $(M_0\cup(-Z),\ga_1)$ to $(M_1,\ga_1)$. Hence there is a morphism
$$F_W:\underline{\rm SHM}(M_0\cup(-Z),\ga_1)\ra\underline{\rm SHM}(M_1,\ga_1).$$
The {\it sutured monopole Floer cobordism map} induced by $\mathcal{W}=(W,Z,\xi)$ is defined as the composition
$$\underline{\rm SHM}(\mathcal{W})=F_W\circ\Phi_{-\xi}:\underline{\rm SHM}(M_0,\ga_0)\ra\underline{\rm SHM}(M_1,\ga_1).$$
\edefn

There are some basic properties of the cobordism map:
\bprop\label{prop_trivial_cobordism}
Suppose Suppose $(M_0,\ga_0)$ is a balanced sutured manifold and $\mathcal{W}=(W,Z,\xi)$ is a suture cobordism from $(M_0,\ga_0)$ to itself so that $W=M_0\times [0,1]$ with $Z=\partial{M}_0\times [0,1]$ and $\xi$ is $I$-invariant. Then we have
$$\underline{\rm SHM}(\mathcal{W})=id:\underline{\rm SHM}(M_0,\ga_0)\ra\underline{\rm SHM}(M_0,\ga_0)$$
\eprop
\bpf
Note the map $F_{W}$ is induced by a cobordism $\widehat{W}$ as in the proof of proposition \ref{prop_map_from_special_cobordism}. In the above settings, however, $\widehat{W}$ is actually diffeomorphic to a product cobordism. Hence the proposition follows from proposition \ref{prop_trivial_gluing_map}.
\epf

\bprop\label{prop_functoriality_of_cobordism_map}
Suppose $(M_0,\ga_0)$, $(M_1,\ga_1)$ and $(M_2,\ga_2)$ are three balanced sutured manifolds. Suppose $\mathcal{W}=(W,Z,\xi)$ is a suture cobordism from $(M_0,\ga_0)$ to $(M_1,\ga_1)$ and $\mathcal{W}'=(W',Z',\xi')$ is a suture cobordism from $(M_1,\ga_1)$ to $(M_2,\ga_2)$. The composition of $\mathcal{W}$ and $\mathcal{W}'$ is a suture cobordism 
$$\mathcal{W}''=(W''=W\cup W',Z''=Z\cup Z',\xi''=\xi\cup\xi')$$
from $(M_0,\ga_0)$ to $(M_2,\ga_2)$. Then we have the equality
$$\underline{\rm SHM}(\mathcal{W}'')=\underline{\rm SHM}(\mathcal{W}')\circ\underline{\rm SHM}(\mathcal{W}):\underline{\rm SHM}(M_0,\ga_0)\ra\underline{\rm SHM}(M_2,\ga_2).$$
\eprop

\bpf
We will not go into details. Suppose we have marked closures
$$\mathcal{D}_0=(Y_0,R_0,r_0,m_0,\eta_0),~\mathcal{D}_1=(Y_1,R_1,r_1,m_1,\eta_1),~\mathcal{D}_2=(Y_2,R_2,r_2,m_2,\eta_2)$$
for $(M_0,\ga_0)$, $(M_1,\ga_1)$ and $(M_2,\ga_2)$ respectively, and the map $\underline{\rm SHM}(\mathcal{W})$ is induced by a cobordism obtained by attaching $4$-dimensional handles $h_1,...,h_n$ and $\tilde{h}_1,...,\tilde{h}_m$ to $Y_0\times[0,1]$ at $Y_0\times \{1\}$. Here $h_1,...,h_n$ correspond to the gluing map $\Phi_{-\xi}$ and $\tilde{h}_1,...,\tilde{h}_m$ correspond to the cobordism map $F_{W}$. Suppose similarly for $\underline{\rm SHM}(\mathcal{W}')$ we have handles $h_1',...,h_s'$ corresponding to the gluing map $\Phi_{-\xi'}$ and $\tilde{h}'_1,...,\tilde{h}'_t$ correspond to the cobordism map $F_{W'}$. Then the composition $\underline{\rm SHM}(\mathcal{W}')\circ\underline{\rm SHM}(\mathcal{W})$ is induced by attaching four sets of $4$-dimensional handles  $h_1,...,h_n$, $\tilde{h}_1,...,\tilde{h}_m$, $h_1',...,h_s'$, $\tilde{h}'_1,...,\tilde{h}'_t$ to $Y_0\times [0,1]$ at $Y_0\times \{1\}$ in the order we wrote them down. Note the attachment of two sets of handles $\tilde{h}_1,...,\tilde{h}_m$ and $h_1',...,h_s'$ can commute with each other because $\tilde{h}_1,...,\tilde{h}_m$ corresponds to handles attached to ${\rm int}(m_1(M_1))\subset Y_1$, while $h_1',...,h_s'$ are attached to $Y_1$ near $m_1(\partial{M})\subset Y_1$ so the two sets of handles are attached disjoint from each other. Then the handles $h_1,...,h_n$ and $h_1',...,h_s'$ are attached first and correspond to the map $\Phi_{-\xi''}$ as in the proof of proposition \ref{prop_functoriality_of_gluing_map}. The handles $\tilde{h}_1,...,\tilde{h}_m$ and $\tilde{h}'_1,...,\tilde{h}'_t$ are attached secondly and correspond to the cobordism map $F_{W''}$ as in the proposition \ref{prop_map_from_special_cobordism}. Hence we get the desired equality:
$$\underline{\rm SHM}(\mathcal{W}'')=\underline{\rm SHM}(\mathcal{W}')\circ\underline{\rm SHM}(\mathcal{W}).$$
\epf

\brem
Intuitively, the three types of maps: cobordism maps, gluing maps and canonical maps all commute with other types. The reason is that for suitable marked closure $\mathcal{D}=(Y,R,r,m,\eta)$, cobordism maps correspond to handles attached in $m({\rm int}(M))\subset Y$, gluing maps correspond to handles attached near $m(\partial{M})\subset Y$ and canonical maps correspond to handles attached in ${\rm int}({\rm im}(r))\subset Y$, and the three regions in $Y$ are pair-wise disjoint.
\erem

\subsection{Duality and turning cobordism around}
Suppose $\mathcal{W}=(W,Z,[\xi])$ is a sutured cobordism from a balanced sutured manifold $(M_1,\ga_1)$ to another $(M_2,\ga_2)$. We can turn the cobordism around, to make another cobordism $\mathcal{W}^{\vee}=(W,Z,[\xi])$ from $(-M_2,\ga_2)$ to $(-M_1,\ga_1)$. Suppose for for $i=1,2$, $\mathcal{D}_i=(Y_i,R_i,r_i,m_i,\eta_i)$ is a marked closure of $(M_i,\ga_i)$, then $\mathcal{D}^{\vee}_i=(-Y_i,-R_i,r_i,-m_i,-\eta_i)$ is a marked closure of $(-M_i,\ga_i)$. Note for a fixed spin${}^c$ structure $\mathfrak{s}$ and smooth $1$-cycle $\eta$ we have a well defined pairing

\begin{equation}\label{eq_pairing_between_closed_manifolds}
	\lgl\cdot,\cdot\rgl: \widehat{HM}(Y,\mathfrak{s};\Gamma_{\eta})\times \widecheck{HM}(-Y,\mathfrak{s};\Gamma_{-\eta})\ra\mathcal{R}.
\end{equation}

Since in sutured monopoles all he spin${}^c$ structures are non-torsion, the pairing \ref{eq_pairing_between_closed_manifolds} induces a pairing
$$\lgl\cdot,\cdot\rgl: SHM(\mathcal{D})\times SHM(\mathcal{D}^{\vee})\ra\mathcal{R}$$
When passing to the projective transitive system and deal with canonical groups or models, the pairing above is well defined up to a unit as we will prove as follows.

\blem\label{lem_dual_for_cobordism_map_between_closed_manifold}
Suppose $W$ is a cobordism from $Y$ to $Y'$ then we can view $W$ as another cobordism $W^{\vee}$ from $-Y'$ to $-Y$. Then the two maps $\widehat{HM}(W)$ and $\widecheck{HM}(W^{\vee})$ are dual to each other with respect to the pairing in (\ref{eq_pairing_between_closed_manifolds}).
\elem

\blem\label{lem_well_definedness_of_pairing_for_sutured_manifold}
Suppose $(M,\ga)$ is a balanced sutured manifold, then there is a pairing well defined up to multiplication by a unit:
\begin{equation}\label{eq_pairing_between_canonical_modules}
	\lgl\cdot,\cdot\rgl:\underline{\rm SHM}(M,\ga)\times \underline{\rm SHM}(-M,\ga)\ra \mathcal{R}.
\end{equation}
\elem
\bpf
First suppose $\mathcal{D}=(Y,R,r,m,\eta)$ and $\mathcal{D}'=(Y',R',r',m',\eta')$ are two marked closures of $(M,\ga)$ of the same genus. Suppose $a\in SHM(\mathcal{D})$ and $b\in SHM(\mathcal{D}^{\vee})$ are two elements in the corresponding homology modules, then we must show that
\begin{equation}\label{eq_well_definedness_of_pairing_for_same_genus}
	\lgl a,b\rgl=\lgl \Phi^g_{\mathcal{D},\mathcal{D}'}(a), \Phi^g_{\mathcal{D}^{\vee},\mathcal{D}'^{\vee}}(b)\rgl,
\end{equation}

where the pairing is the one in (\ref{eq_pairing_between_closed_manifolds}).

We prove here only the case when there is a curve $\al\subset R$ so that after doing a $(+1)$ surgery along $r(\al\times\{0\})\subset Y$ with respect to the $r(R\times\{0\})$-framing, we get a manifold diffeomorphic to $Y'$. The general case will follow from a similar argument and the functoriality of the canonical map $\Phi^g$.

Suppose there is a curve $\al'\subset R$ parallel to $\al$ but is disjoint from $\al$. Since the Dehn surgery is supported in arbitrary small neighborhood of $\al$, we can assume that $r(\al'\times\{0\})\subset Y'$.

In this case, there is a cobordism $W^+$ from $Y$ to $Y'$ obtained by attaching a $4$-dimensional $2$-handle, with $(+1)$-framing with respect to the $r(R\times\{0\})$-surface framing, to $Y\times[0,1]$ along $r(\al\times\{0\})\times\{1\}\subset Y\times\{1\}$, and
$$\Phi^g_{\mathcal{D},\mathcal{D'}}=HM(W^+).$$

On $\mathcal{D}^{\vee}$, the surgery is still a $(+1)$-surgery, but we  so there is a cobordism $W^-$ from $-Y'$ to $-Y$ obtained by gluing a $4$-dimensional $2$-handle, with $(+1)$-framing with respect to $-r'(R'\times\{0\})$, to the curve $r(\al'\times\{0\})\times\{1\}\subset Y'\times\{1\}$, and
$$\Phi^g_{\mathcal{D}^{\vee},\mathcal{D'}^{\vee}}=HM(W^-)^{-1}.$$

We actually have that $W^+$ and $W^-$ are diffeomorphic by an orientation preserving diffeomorphism, so $W^-$ can be viewed as turning $W^+$ around. As a result, by lemma \ref{lem_dual_for_cobordism_map_between_closed_manifold}, we have
\beq
\lgl \Phi^g_{\mathcal{D},\mathcal{D}'}(a), \Phi^g_{\mathcal{D}^{\vee},\mathcal{D}'^{\vee}}(b)\rgl&=\lgl HM(W^+)(a),HM(W^-)^{-1}(b)\rgl\\
&=\lgl a,HM(W^-)\circ HM(W^-)^{-1}(b)\rgl\\
&=\lgl a,b\rgl.
\eeq
Hence (\ref{eq_well_definedness_of_pairing_for_same_genus}) is proved.

Now suppose $\mathcal{D}=(Y,R,r,m,\eta)$ and $\mathcal{D}'=(Y',R',r',m',\eta')$ are two marked closures for $(M,\ga)$ so that $g(\mathcal{D}')=g(\mathcal{D})+1$. Then we need to show that

\begin{equation}\label{eq_well_definedness_of_pairing_for_different_genus}
	\lgl a,b\rgl=\lgl \Phi^{g,g+1}_{\mathcal{D},\mathcal{D}'}(a), \Phi^{g,g+1}_{\mathcal{D}^{\vee},\mathcal{D}'^{\vee}}(b)\rgl,
\end{equation}

Since we have dealt with the case of the same genus, we can discuss only the special case as follows: there are two disjoint oriented embedded tori $T_1,T_2\subset Y'$ so that

(1). For $i=1,2$, $T_i\cap m'(M)=\emptyset$.

(2). For $i=1,2$, $T_i\cap r'(R'\times[-1,1])=r'(c_i\times[-1,1])$ where $c_i\subset R'$ is an embedded oriented circle, and the two circle $c_1$ and $c_2$ together cut $R'$ into two oriented parts $R_1'$ and $R_2'$, so that 
$$c_1\cup c_2=\partial{R}'_1=-\partial{R}'_2~{\rm and}~ R'_2\cong \Sigma_{1,2},$$ where $\Sigma_{1,2}$ is the compact oriented surface of genus $1$ and having two boundary component.

(3). $T_1$ and $T_2$ cut $Y'$ into two parts $Y'_1$ and $Y'_2$ so that
$$T_1\cup T_2=\partial{Y}_1=-\partial{Y}'_2~{\rm and}~m'(M)\subset Y_1'.$$

(4). For $i=1,2$, $\eta'$ intersects $R'_i$ in an oriented, non-boundary-parallel properly embedded arc $\eta'_i$.

Suppose for $i=1,2$, $p_i=c_i\cap\eta_i$ and pick an orientation reversing diffeomorphism $f:c_1\ra c_2$ sending $p_1$ to $p_2$. Choose an orientation reversing diffeomorphism $h:T_1\ra T_2$ so that for $i=1,2$
$$h|_{r'(c_i\times[-1,1])}=(r')^{-1}\circ(f\times id)\circ r'.$$
We can use $h$ to glue the two boundary components of $Y_1'$ to get a closed manifold $Y_1$ and do the same thing for $Y_2'$ to get a closed manifold $Y_2$. As done in \cite{baldwin2015naturality} there is a natural way to get a closure $\mathcal{D}_1=(Y_1,R_1,r_1,m_1,\eta_1)$ for $(M,\ga)$ and a closure $\mathcal{D}_2=(Y_2,R_2,r_2,m_2,\eta_2)$ so that $Y_2$ is a fibration over $S^1$ with fibres diffeomorphic to $R_2$. 

As we have already deal with the case of same genus, we can assume that the two marked closures $\mathcal{D}$ and $\mathcal{D}_1$ are the same. Then we can describe the canonical maps $\Phi^{g,g+1}_{\mathcal{D},\mathcal{D}'}$ and $\Phi^{g,g+1}_{\mathcal{D}^{\vee},\mathcal{D}'^{\vee}}$ as follows. Pick the surface $U$ depicted in figure \ref{fig_excision_cobordism}. Glue the three part $Y_1'$, $T_1\times U$ and $Y_2'$ together using $h$ just as depicted by figure \ref{fig_excision_cobordism}. The result is a cobordism $W^+$ from $Y=Y_1$ disjoint union $Y_2$ to $Y'$. This cobordism will induce the canonical map $\Phi^{g,g+1}_{\mathcal{D},\mathcal{D}'}$. The same cobordism, with the reversed orientation, will be a cobordism $W^-=-W^+$ from $(-Y)\sqcup(-Y_2)$ to $-Y'$ and it will induce the canonical map $\Phi^{g,g+1}_{\mathcal{D}^{\vee},\mathcal{D}'^{\vee}} $. If we turn $W^+$ around, it will become a cobordism $W^{\vee}$ from $-Y'$ to $-Y\sqcup(-Y_2)$ and induce a dual map by lemma \ref{lem_dual_for_cobordism_map_between_closed_manifold}. Then the equality (\ref{eq_well_definedness_of_pairing_for_different_genus}) will follow from the fact that the cobordism $W^{\vee}\cup W^-$ will induce the identity map up to multiplication by a unit, which is proved in \cite{kronheimer2010knots}.
\epf

There is a simpler way to describe the gluing map.

Suppose $(M',\ga')$ is a balanced sutured manifold and $(M,\ga)$ is a sutured submanifold. Suppose $Z=M'\backslash{\rm int}(M)$ and $\xi$ is a contact structure on $Z$ so that $\partial{Z}$ is convex with dividing set $\ga\cup\ga'$. Suppose $Z$ has a contact handle decomposition relative to $M$. That is, there are contact handles $h_1,...,h_n$ so that if we attach them to $(M,\ga)$, then we will get $(M',\ga')$. Suppose $h_1,...,h_m$ are all $0$- and $1$-handles and $h_{m+1},...,h_{n}$ are all $2$- and $3$-handles. Suppose $(M_1,\ga_1)$ is the result of attaching all $h_1,...,h_m$ to $(M,\ga)$. Let $W=M'\times[0,1]$, and let $M_2=\partial{W}\backslash (M_1\times\{0\})$ with suitable orientation. We can view $W$ as a cobordism from $(M_1,\ga_1)$ to $(M_2,\ga_1)$ with sutured surface $(S=\partial{M}_1\times\{0\},\ga_1)$. If we do closing up along $S$, we will get two marked closures $\mathcal{D}_1=(Y_1,R,r,m_1,\eta)$ and $\mathcal{D}_2=(Y_2,R,r,m_2,\eta)$ for $(M_1,\ga_1)$ and $(M_2,\ga_2)$ respectively and a cobordism $\widehat{W}$ from $Y_1$ to $Y_2$. 
\bprop\label{prop_another_description_of_gluing_map}
Under the above settings, the marked closure $\mathcal{D}_1$ is also a marked closure for $(M,\ga)$ so there is a map
$$\Phi:\underline{\rm SHM}(-M,-\ga)\ra \underline{\rm SHM}(-M_1,-\ga_1).$$
The marked closure $\mathcal{D}_2$ is also a marked closure for $(M',\ga')$ so there is a map
$$\Psi:\underline{\rm SHM}(-M',-\ga')\ra \underline{\rm SHM}(-M_2,-\ga_2).$$
The gluing map can be written as
$$\Phi_{\xi}=\Psi^{-1}\circ F_{-W}\circ \Phi.$$
\eprop

\bpf
From proposition \ref{prop_equivalence_between_handle_and_general_gluing} we know that the gluing map is actually equal to
$$\Phi_{\xi}=C_{h_n}\circ...\circ C_{h_1}.$$
Since $(M_1,\ga_1)$ is gotten from $(M,\ga)$ by attaching a few $0$- and $1$-handles, the marked closure $\mathcal{D}$ for $(M_1,\ga_1)$ must also be one for $(M,\ga)$ and hence $\Phi$ is just the composition
$$\Phi=C_{h_m}\circ...\circ C_{h_1}.$$

Let $W_1=M_1\times[0,1]\subset W$ be the product. Still let $S=\partial{M}_1\times\{0\}$ and let $M_3=\partial{W}_1\backslash (M_1\times\{0\})$ with suitable orientation. Then we can view $W$ as a cobordism from $(M_1,\ga_1)$ to $(M_3,\ga_3)$ with sutured surface $(S,\ga_1)$. When doing the same closing up along $S$ as above, we get two marked closures $\mathcal{D}_1$ and $\mathcal{D}_3=(Y_3,R,r,m_3,\eta)$ for $(M_1,\ga_1)$ and $(M_3,\ga_1)$ respectively and a cobordism $\widehat{W}_1$ from $Y_1$ to $Y_3$. If we write
$$h_j=(\phi_j,S_j,D^3_j,\delta_j),$$
then we can see that $\hat{W}$ is gotten from $\widehat{W}_1$ by attaching all $D^3_j\times [0,1]$ to $\widehat{W}_1$ via maps
$$\phi_j\times id: S_j\times[0,1]\ra \partial{M}_1\times[0,1]\subset M_2\subset Y_2\subset \partial{\widehat{W}}_1.$$
This exactly the way we define $2$- and $3$-handle attaching maps in definition \ref{defn_two_handle_gluing_map2} and definition \ref{defn_three_handle_gluing_map2}. Hence we have
$$\Psi^{-1}\circ F_{-W}=C_{n}\circ...\circ C_{m+1}$$
and we are done.
\epf

\bcor\label{cor_another_description_of_cobordism_map}
Suppose $\mathcal{W}=(W,Z,[\xi])$ is a sutured cobordism from $(M_1,\ga_1)$ to $(M_2,\ga_2)$. Suppose $S\subset Z$ is chosen as in proposition \ref{prop_another_description_of_gluing_map} and $\ga_1'$ correspondingly. Suppose $S$ separates $\partial{W}$ into two parts $M_1'$ and $M_2'$, so that $M_i'$ contains $M_i$ and is oriented in the same way as $M_i$. We can view $W$ as a cobordism from $(M_1',\ga_1')$ to $(M_2',\ga_1')$ with sutured surface $(S,\ga_1')$. If we do closing up along $S$, we get two marked closures $\mathcal{D}_1'=(Y_1',R,r,m_1',\eta)$ and $\mathcal{D}_2'=(Y_2',R,r,m_2',\eta)$ for $(M_1',\ga_1')$ and $(M_2',\ga_1')$ respectively. As above we have

$$\Phi:\underline{\rm SHM}(M_1,\ga_1)\ra \underline{\rm SHM}(M_1',\ga_1'),$$
$$\Psi:\underline{\rm SHM}(M_2,\ga_2)\ra \underline{\rm SHM}(M_2',\ga_1').$$
Then we can actually write the cobordism map to be
$$\shm(\mathcal{W})=\Psi^{-1}\circ F_{W}\circ\Phi.$$
\ecor

\bpf
We can decompose the sutured cobordism $\mathcal{W}$ as a union of two: 
$$\mathcal{W}=\mathcal{W}^b\cup\mathcal{W}^s.$$

Here $\mathcal{W}^s$ is a special cobordism whose underline manifold is $W^s=W$ but with sutured surface $(S'=\partial{M}_2,\ga_2)$. The cobordism $\mathcal{W}^b$ is a special cobordism whose underlining manifold is $W^b=(M\cup(-Z))\times[0,1]$ but with sutured surface $(S\times\{0\}\subset (M\cup(-Z))\times[0,1], \ga_1')$. The cobordism $W$ can be viewed as a union $W^s\cup W^b$ with sutured surface $(S,\ga_1')$. From proposition \ref{prop_map_from_special_cobordism} we know that $W^s$ can be viewed as gotten from $(M_1\cup(-Z))\times[0,1]$ by attaching some $4$-dimensional handles $h^4_1,...,h^4_l$ to ${\rm int}(M_1\cup(-Z))\times\{1\}$.  Hence the result $\widehat{W}$ of doing closing up along $S$ for $W$ is the same as doing closing up along $S$ for for $W^s\cup W^b$. The result of the later can be described as follows. When doing closing up along $S$ for $W^b$, we get two marked closure $\mathcal{D}_1=(Y_1,R,r,m_1,\eta)$ and $\mathcal{D}=(Y,R,r,m,\eta)$ for $(M_1,\ga_1)$ and $(M_1\cup(-Z),\ga_2)$ respectively, and a cobordism $\hat{W}^b$ from $Y_1$ to $Y$. Now $\widehat{W}^b$ can be thought of as obtained from $Y_1\times [0,1]$ by attaching some $4$-dimensional $2$- and $3$- handles which correspond to the gluing map $\Phi_{-\xi}$. When adding the contribution from $W^s$, we know that $\widehat{W}$ can be viewed as obtained from $\widehat{W}^b$ by attaching $h^4_1,...,h^4_l$ to $Y\subset\partial{\hat{W}^b}$. This description also exists in the construction of cobordism map. Hence we know that $\widehat{W}$ indeed induces the cobordism map $\shm({\mathcal{W}})$.

\epf

Now we can describe the relation between $\shm(\mathcal{W})$ and $\shm(\mathcal{W}^{\vee})$ as follows:

\bcor\label{cor_turning_sutured_cobordism_around}
Suppose $\mathcal{W}=(M,Z,[\xi])$ is a sutured cobordism from $(M_1,\ga_1)$ to $(M_2,\ga_2)$. The same cobordism can be also viewed as a cobordism $\mathcal{W}^{\vee}$ from $(-M_2,\ga_2)$ to $(-M_1,\ga_1)$. Then the cobordism map $\shm(\mathcal{W})$ and $\shm(\mathcal{W}^{\vee})$ are dual with respect to the pairing (\ref{eq_pairing_between_canonical_modules}).
\ecor
\bpf
If put the cobordism $W$ up-side-down, then the distinguising surface $S$ is unchanged (but in $Z$, a $i$-handle becomes a $(3-i)$-handle). Hence from corollary \ref{cor_another_description_of_cobordism_map}, $F_{\mathcal{W}}$ is induced by a cobordism $\hat{W}$ while $F_{\mathcal{W}^{\vee}}$ is induced by a cobordism $\hat{W}^{\vee}$ which is obtained by putting $\hat{W}$ up-side-down. Hence the conclusion follows from lemma \ref{lem_dual_for_cobordism_map_between_closed_manifold}.
\epf

There is a question related to the trace and co-trace cobordism. Suppose $(M,\ga)$ is a balanced sutured manifold and $\mathcal{W}=(W=M\times[0,1],Z=\partial{M}\times[0,1],[\xi])$ is the sutured cobordism from $(M\sqcup(-M),\ga\cup(-\ga))$. Here $\xi$ is a $[0,1]$-invariant contact structure on $Z$ so that $\partial{M}$ is convex with respect to $\frac{\partial}{\partial{t}}$ and $\ga$ is the corresponding dividing set. Let $\mathcal{R}$ be the ring with which we build the local coefficient, then we would like to ask the following question:

\begin{quest}\label{quest_trace_and_cotrace}
	How to describe the cobordism map 
	$$\shm(\mathcal{W}): \underline{\rm SHM}(M\sqcup (-M),\ga\cup \ga)\ra \mathcal{R}?$$
\end{quest}

Note from K\"unneth formula, there is a map
$$i:\underline{\rm SHM}(M\sqcup (-M),\ga\cup \ga)\ra\underline{\rm SHM}(M,\ga)\otimes\underline{\rm SHM}(-M,\ga).$$
Also there is a canonical map
$$\mathop{tr}:\underline{\rm SHM}(M,\ga)\otimes\underline{\rm SHM}(-M,\ga)\ra \mathcal{R}$$
defined as 
$$\mathop{tr}(a\otimes b)=b(a),$$
since $\underline{\rm SHM}(-M,\ga)$ is the dual of $\underline{\rm SHM}(M,\ga)$. We make the following conjecture:

\begin{conj}\label{conj_trace_cobordism}
With the above settings, we have
\begin{equation}\label{eq_conjecture_for_trace_cobordism_map}
\shm(\mathcal{W})=\mathop{tr}\circ i.	
\end{equation}
\end{conj}

\section{A brief discussion on Instanton}
The constructions in section 3-5 can be applied to instanton sutured manifolds. The following definition is from \cite{baldwin2016instanton}.

\bdefn\label{defn_marked_odd_closure}
Suppose $(M,\ga)$ is a balanced sutured manifold, then a {\it marked odd closure} of $(M,\ga)$ is a sextuple $\mathcal{D}=(Y,R,m,r,\eta,\al)$ so that

(1). The quintuple $(Y,R,m,r,\eta)$ is a marked closure of $(M,\ga)$ defined as in definition \ref{defn_marked_closure}

(2). We have $\al$ being a curve disjoint from ${\rm im}(m)$ and intersects $r(R\times[-1,1])$ in the form $r(\{p\}\times[-1,1]$ for some point $p\in R$.
\edefn

Now suppose $\mathcal{D}=(Y,R,m,r,\eta,\al)$ is a marked odd closure of a balanced sutured manifold $(M,\ga)$ we can pick a Hermitian line bundle $\omega$ over $Y$ such that $c_1(\omega)$ is dual to the curve $\al\cup\eta$. Let $E$ be a $U(2)$-bundle over $Y$ with a bundle isomorphism $\rho: \Lambda^2E\ra L$. With such data we could define instanton Floer homology $I_{*}(Y)_{\omega}$ on $Y$. Follow from the definition in \cite{kronheimer2010knots}, we can define
$$SHI(\mathcal{D})=I_{*}(Y|r(R\times\{0\}))_{\omega},$$
where $I_{*}(Y|r(R\times\{0\}))_{\omega}$ means the generalized eigenspace of $\mu(r(R\times\{0\}))$ in $I_{*}(Y)_{\omega}$ with eigenvalue $2g(R)-2$. In \cite{baldwin2015naturality}, Baldwin and Sivek construct canonical maps between marked odd closures and the sutured instanton Floer homology becomes a projective transitive system of $\mathbb{C}$-modules. In \cite{baldwin2016instanton}, they also construct contact handle gluing maps for instanton Floer homology and the construction in this paper would be applied to instanton and we have:
\bthm
For sutured instanton Floer homology, we have:

(1). The handle gluing maps constructed by Baldwin and Sivek satisfy similar cancelation and invaraint properties as in lemmas \ref{lem_zero_one_cancelation}, \ref{lem_one_two_handle_cancelation}, \ref{lem_two_three_handle_cancelation}, \ref{lem_independent_of_order}, \ref{lem_commute_with_diffeomorphisms}, \ref{lem_gluing_disjoint_union}.

(2). There are well defined (up to multiply by a non-zero complex number) gluing maps for sutured instanton Floer homology and it satisfies similar properties as in propositions \ref{prop_functoriality_of_gluing_map}, \ref{prop_trivial_gluing_map}, \ref{prop_equivalence_between_handle_and_general_gluing}.

(3). There are well defined (up to multiply by a non-zero complex number) gluing maps for sutured instanton Floer homology and it satisfies similar properties as in propositions \ref{prop_functoriality_of_cobordism_map}, \ref{prop_trivial_cobordism}.
\ethm

\brem
This will give a confirmative answer to conjecture 1.8 in \cite{baldwin2016instanton} where Baldwin and Sivek conjectures that the gluing maps is independent of the handle decomposition.
\erem

At last we want to give an alternative definition of the contact invariant defined in \cite{baldwin2016instanton}. The original definition used partial open book decompositions of contact balanced sutured manifolds. Yet partial open book decompositions only involve $0$-, $1$- and $2$- handles so it is only expected but not proved that the contact element is also preserved by attaching a contact $3$-handle. Now we can use the gluing map to define the contact element:

\bdefn\label{defn_contact_element_in_shi}
Suppose $(M,\ga)$ is a balanced sutured manifold and $\xi$ is a contact structure on $M$ so that $\partial{M}$ is a convex surface and $\ga$ is the dividing set. We can define the {\it contact element} of $(M,\ga,\xi)$ as follows. Suppose $D\subset {\rm int}(M)$ is a Darboux ball in $M$, let $\delta\subset \partial{D}$ be the dividing set on $\partial{D}$. Let $Z=M\backslash{\rm int}(D)$, we have a gluing map
$$\Phi_{\xi}:\underline{\rm SHI}(-D,-\delta)\ra \underline{\rm SHI}(-M,-\ga).$$
Then the contact element $\phi(M,\ga,\xi)\in \underline{\rm SHI}(-M,-\ga)$ is defined as
$$\phi(M,\ga,\xi)=\Phi_{\xi}(1),$$
where $1\in \underline{\rm SHI}(-D,-\delta)$ is a generator of the canonical module.
\edefn

\bprop
Suppose $(M,\ga)$ is a balanced sutured manifold and $\xi$ is a contact structure on $M$ so that $\partial{M}$ is a convex surface with dividing set $\ga$. Then

(1). The contact invariant defined as in definition \ref{defn_contact_element_in_shi} is equivalent to the contact element defined by Baldwin and Sivek in \cite{baldwin2016contact}.

(2). The contact element is preserved under the gluing map $\Phi_{\xi}$.
\eprop

As discussed in \cite{kronheimer2010knots}, if we have a closed $3$-manifold $Y$ and we dig a $3$-ball to create a spherical boundary with one simple closed curve as sutures, then the sutured instanton homology can be identified with the instanton Floer homology of a suitable admissible bundle over $Y\sharp T^3$. So Baldwin and Sivek's construction would result in a contact element for the closed $3$-manifold $Y\sharp T^3$.

\begin{quest}
Can we re-construct the contact element in the classical instanton Floer homology theory? Would this element be preserved by exact symplectic cobordism as so in the monopole settings?
\end{quest}

\bibliography{Index}
\end{document}